\documentclass{article}

\usepackage{placeins}
\usepackage[preprint]{nips_2019}
\usepackage{graphicx,verbatim,color}
\usepackage{enumitem}
\usepackage{natbib}
\usepackage{amsmath}
\usepackage{amssymb}
 \usepackage{multirow}
 \usepackage{mathrsfs}
\usepackage{amsfonts}
\usepackage{algorithm}
\usepackage{algorithmic}
\usepackage{graphicx, subfig}
\usepackage{epsfig}
\usepackage{epstopdf}
\usepackage{amsthm}\usepackage{dsfont}\usepackage{array}\usepackage{mathrsfs}\usepackage{comment}
\usepackage{hyperref}
\usepackage{breakurl}
\usepackage[utf8]{inputenc} 
\usepackage[T1]{fontenc}    
\usepackage{hyperref}       
\usepackage{url}            
\usepackage{booktabs}       
\usepackage{amsfonts}       
\usepackage{nicefrac}       
\usepackage{microtype}      
\newcommand{\eq}[1]{Eq.~(\ref{eq:#1})}
\newcommand{\tabincell}[2]{\begin{tabular}{@{}#1@{}}#2\end{tabular}}

\newcommand{\note}[1]{{\textbf{\color{red}#1}}}

\newcommand{\trace}{{\rm {trace}}}
\newcommand{\rank}{{\rm {rank}}}

\newcommand{\mR}{\mathbb{R}}

\newcommand{\half}{\frac{1}{2}}

\newcommand{\lc}{\left(}
\newcommand{\rc}{\right)}

\newcommand{\lbig}{\left\{}
\newcommand{\rbig}{\right\}}

\newcommand{\lb}{\lambda}

\newcommand{\nb}{\nabla}
\newcommand{\ep}{\epsilon}

\newcommand{\bt}{\beta}
\newcommand{\dt}{\delta}
\newcommand{\ka}{\kappa}
\newcommand{\gm}{\gamma}
\newcommand{\Gm}{\Gamma}
\newcommand{\Dt}{\Delta}

\newcommand{\minl}{\mathop{\min}\limits_}

\newcommand{\ra}{\rightarrow}

\newcommand{\la}{\leftarrow}

\newcommand{\taux}{x_\tau}

\newcommand{\diag}{{\rm diag}}

\newcommand{\mG}{\mathscr{G}}
\newcommand{\mF}{\mathscr{F}}
\newcommand{\mathS}{\mathscr{S}}
\newcommand{\mathI}{\mathfrak{I}}
\newcommand{\mathT}{\mathscr{T}}

\newcommand{\mT}{\mathcal{T}}
\newcommand{\mM}{\mathcal{M}}

\newcommand{\mP}{\mathcal{P}}

\newcommand{\grad}{{\rm grad}}
\newcommand{\Exp}{{\rm Exp}}
\newcommand{\Expn}{{\rm Exp}^{-1}}

\newcommand{\ys}[1]{\noindent{\textcolor{blue}{\{{\bf YS:} #1\}}}}

\newcommand{\wc}{\widecheck}

\DeclareFontFamily{U}{mathx}{\hyphenchar\font45}
\DeclareFontShape{U}{mathx}{m}{n}{
      <5> <6> <7> <8> <9> <10>
      <10.95> <12> <14.4> <17.28> <20.74> <24.88>
      mathx10
      }{}
\DeclareSymbolFont{mathx}{U}{mathx}{m}{n}
\DeclareMathAccent{\widecheck}{\mathalpha}{mathx}{"71}

\title{Escaping from saddle points on Riemannian manifolds}
\author{Yue Sun\\
University of Washington\\Seattle, WA, 98105\\ \texttt{yuesun@uw.edu}\And
Nicolas Flammarion\\
University of California, Berkeley\\Berkeley, CA, 94720\\ \texttt{flammarion@berkeley.edu}\And
Maryam Fazel\\
University of Washington\\Seattle, WA, 98105\\ \texttt{mfazel@uw.edu}}

\begin{document}

\theoremstyle{plain}\newtheorem{assumption}{Assumption}\theoremstyle{plain}\newtheorem{lemma}{\textbf{Lemma}}\newtheorem{condition}{\textbf{Condition}}\newtheorem{theorem}{\textbf{Theorem}}\newtheorem{corollary}{\textbf{Corollary}}\newtheorem{example}{\textbf{Example}}\newtheorem{definition}{\textbf{Definition}}\newtheorem{conjecture}{\textbf{Conjecture}}
\newtheorem{claim}{\textbf{Claim}}\newtheorem{remark}{\textbf{Remark}}\newtheorem{question}{\textbf{Question}}
\theoremstyle{definition}

\maketitle

 \begin{abstract}
We consider minimizing a nonconvex, smooth function $f$ on a Riemannian manifold $\mathcal{M}$. We show that a perturbed version of Riemannian gradient descent algorithm converges to a second-order stationary point (and hence is able to \emph{escape} saddle points on the manifold). The rate of convergence depends as $1/\epsilon^2$ on the accuracy $\epsilon$, which matches a rate known only for unconstrained smooth minimization. 
The convergence rate depends  polylogarithmically on the manifold dimension $d$, hence is almost dimension-free. 
The rate also has a polynomial dependence on the parameters describing the curvature of the manifold and the smoothness of the function. 
While the unconstrained problem (Euclidean setting) is well-studied, our result is the first to prove such a rate for nonconvex, manifold-constrained problems.
\end{abstract}
\vspace{-0.4cm}
\section{Introduction}
\vspace{-0.2cm}
We consider minimizing a non-convex smooth function on a smooth manifold $\mM$, 
\begin{align}\label{eq:problem}
\mbox{minimize}_{x\in\mM}~~ f(x),
\end{align}
where $\mM$ is a $d$-dimensional smooth manifold\footnote{Here $d$ is the dimension of the manifold itself; we do not consider $\mM$ as a submanifold of a higher dimensional space. For instance, if $\mM$ is a 2-dimensional sphere embedded in $\mathbb{R}^3$, its dimension is $d=2$.}, and $f$ is twice differentiable, with a Hessian that is $\rho$-Lipschitz (assumptions are formalized in section 4). 
This framework includes a wide range of fundamental problems (often non-convex), such as PCA \citep{edelman1998geometry}, dictionary learning \citep{sun2017complete}, low rank matrix completion \citep{boumal2011rtrmc}, and tensor factorization \citep{ishteva2011best}. Finding the global minimum to \eq{problem} is in general NP-hard; our goal is to find an approximate second order stationary point with first order optimization methods. We are interested in first-order methods because they are extremely prevalent in machine learning, partly because computing Hessians is often too costly. It is then important to understand how first-order methods fare when applied to nonconvex problems, and there has been a wave of recent interest on this topic since 
\citep{Ge}, as reviewed below. 


In the Euclidean space,
it is known that with random initialization, gradient descent avoids saddle points asymptotically \citep{nonconvergence, OnlyMin}. \citet{lee2017first} (section 5.5) show that this is also true on smooth manifolds,
although the result is expressed in terms of nonstandard manifold smoothness measures. 
Also, importantly, this line of work does not give quantitative rates for the algorithm's behaviour near saddle points.

\citet{ExpTime} show gradient descent can be \emph{exponentially slow} in the presence of saddle points. To alleviate this phenomenon,
it is shown that for a $\bt$-gradient Lipschitz, $\rho$-Hessian Lipschitz function, cubic regularization \citep{CubicR} and perturbed gradient descent \citep{Ge,GeSham} converges to $(\ep, -\sqrt{\rho\ep})$ local minimum \footnote{defined as $x$ satisfying $\|\nb f(x)\|\le \ep$, $\lb_{\min} \nb^2 f(x)\ge -\sqrt{\rho\ep}$} in polynomial time, and momentum based method accelerates \citep{HeavyBall}. 
Much less is known about  
inequality constraints:
\citet{Lee_constrained} and \citet{mokhtari_constrained} discuss second order convergence for general inequality-constrained problems, where they need an
NP-hard subproblem (checking the co-positivity of a matrix) to admit a polynomial time approximation algorithm. 
However such an approximation exists only under very restrictive assumptions. 

An orthogonal line of work is optimization on Riemannian manifolds. \citet{OptMan} provide comprehensive background, showing how algorithms such as gradient descent, Newton and trust region methods can  be implemented on Riemannian manifolds, together with asymptotic convergence guarantees to first order stationary points. \citet{Sra_geo_convex} provide global convergence guarantees for first order methods when optimizing geodesically convex functions. \citet{bonnabel2013stochastic} obtains the first asymptotic convergence result for stochastic gradient descent in this setting, which is further extended by \citet{manifold_first,svrg_manifold,khuzani2017stochastic}. If the problem is non-convex, or the Riemannian Hessian is not positive definite, one can use second order methods to escape from saddle points. \citet{Rie_trust_region} shows that Riemannian trust region method converges to a second order stationary point in polynomial time~\citep[see, also,][]{NIPS2018_7679,MR3826674,zhang2018cubic}. But this method requires a Hessian oracle, whose complexity is $d$ times more than computing gradient. In Euclidean space, trust region subproblem can be sometimes solved via a Hessian-vector product oracle, whose complexity is about the same as computing gradients. \citet{agarwal2018adaptive} discuss its implementation on Riemannian manifolds, but not clear about the complexity and sensitivity of Hessian vector product oracle on manifold. 

The study of the convergence of gradient descent for non-convex Riemannian problems is previously done only in the Euclidean space by modeling the manifold with equality constraints. \citet[][ Appendix B]{Ge} prove that stochastic projected gradient descent methods converge to  second order stationary points in polynomial time (here the analysis is not geometric, and depends on the algebraic representation of the equality constraints).  \citet{Yue_projected} proves perturbed projected gradient descent converges with a comparable rate to the unconstrained setting~\citep{GeSham} (polylog in dimension). 
The paper applies projections from the ambient Euclidean space to the manifold and analyzes the iterations under the Euclidean metric. This approach loses the geometric perspective enabled by Riemannian optimization, and cannot explain convergence rates in terms of inherent quantities such as the sectional curvature of the manifold.

After finishing this work, we found the recent and independent paper
\citet{Criscitiello2019Efficiently} which gives a similar convergence analysis result for a related  perturbed Riemannian gradient method. We point out a few differences:
(1) In \citet{Criscitiello2019Efficiently} Lipschitz assumptions are made on the pullback map $f\circ \mathrm{Retr}$. While this makes the analysis simpler, it lumps the properties of the function and the manifold together, and the role of the manifold's curvature is not explicit. In contrast, our rates are expressed in terms of the function's smoothness parameters and the sectional curvature of the manifold separately, capturing the geometry more clearly. 
(2) The algorithm in \citet{Criscitiello2019Efficiently} uses two types of iterates (some on the manifold but some taken on a tangent space), whereas all our algorithm steps are directly on the manifold, which is more natural.
(3) To connect our iterations with intrinsic parameters of the manifold, we use the exponential map instead of the more general retraction used in \citet{Criscitiello2019Efficiently}. 

{\bf Contributions.} We provide convergence guarantees for perturbed first order Riemannian optimization methods to second-order stationary points (local minimum). We prove that as long as the function is appropriately smooth and the manifold has bounded sectional curvature, a perturbed Riemannian gradient descent algorithm escapes (an approximate) saddle points with a rate of $1/\epsilon^2$, a polylog dependence on the dimension of the manifold (hence almost dimension-free), and a polynomial dependence on the smoothness and curvature parameters. This is the first result showing such a rate for Riemannian optimization, and the first to relate the rate to geometric parameters of the manifold.

Despite analogies with the unconstrained (Euclidean) analysis and with the Riemannian optimization literature, 
the technical challenge in our proof goes beyond combining two lines of work:
we need to  
analyze the interaction between the first-order method and the second order structure of the manifold to obtain second-order convergence guarantees that depend on the manifold curvature. Unlike in Euclidean space, the curvature affects the Taylor approximation of gradient steps. On the other hand, unlike in the local rate analysis in first-order Riemannian optimization, our second-order analysis requires more refined properties of the manifold structure (whereas in prior work, first order oracle makes enough progress for a local convergence rate proof, see Lemma \ref{lem:first_order_progress}), and second order algorithms such as \citep{Rie_trust_region} use second order oracles (Hessian evaluation). See section 4 for further discussion.
\vspace{-0.2cm}
\section{Notation and Background}
\vspace{-0.2cm}
We consider a complete\footnote{Since our results are local, completeness is not necessary and our results can be easily generalized, with extra assumptions on the injectivity radius.}, smooth, $d$ dimensional Riemannian manifold $(\mM,\mathfrak g)$, equipped with a Riemannian metric $\mathfrak g$, 
and we denote by $\mT_x\mM$ its tangent space at $x\in \mM$ (which is a vector space of dimension $d$). We also denote by $\mathbb{B}_x(r)=\{v\in\mT_x\mM,\Vert v \Vert\le r\}$ the ball of radius $r$ in $\mT_x\mM$ centered at $0$. 
At any point $x\in \mM$, the metric $\mathfrak g$ induces a natural inner product on the tangent space denoted by $\langle\cdot,\cdot\rangle: \mT_x\mM\times \mT_x\mM\ra \mathbb{R}$. We also consider  the Levi-Civita connection  $\nb$ (\citealp{OptMan}, Theorem 5.3.1).  The Riemannian curvature tensor is denoted by $R(x)[u,v]$ where $x\in\mM$, $u,v\in\mT_x\mM$ and is defined in terms of the connection $\nb$~\citep[][Theorem 5.3.1]{OptMan}. The sectional curvature $K(x)[u,v]$ for $x\in\mM$ and $u,v\in\mT_x\mM$ is then defined in 
\citet[][Prop. 8.8]{LeeJohn97}.
\begin{align*}
    K(x)[u,v] = \frac{\langle R(x)[u,v]u,v\rangle}{\langle u,u\rangle\langle v,v\rangle-\langle u,v\rangle^2},\ x\in\mM, \ u,v\in\mT_x\mM.
\end{align*}


Denote the distance (induced by the Riemannian metric) between two points in $\mM$ by $d(x,y)$. A geodesic $\gm:\mR\ra\mM$ is a constant speed curve whose length is equal to $d(x,y)$, so it is the shortest path on manifold linking $x$ and $y$. $\gm_{x\ra y}$ denotes the geodesic from $x$ to $y$ (thus $\gm_{x\ra y}(0) = x$ and $\gm_{x\ra y}(1)=y$).


The exponential map $\Exp_x(v)$ maps $v\in\mT_x\mM$ to $y\in \mM$ such that there exists a geodesic $\gamma$ with $\gamma(0)=x$, $\gamma(1)=y$ and $\frac{d}{dt}\gamma(0)=v$. The injectivity radius at point $x\in\mM$ is the maximal radius $r$ for which the exponential map is a diffeomorphism on $\mathbb{B}_x(r)\subset\mT_x\mM$. The injectivity radius of the manifold, denoted by $\mathI$, is the infimum of the injectivity radii at all points. Since the manifold is complete, we have $\mathI>0$. When $x,y\in\mM$ satisfies $d(x,y)\le \mathI$, 
the exponential map admits an inverse $\Exp_x^{-1}(y)$, which satisfies $d(x,y)=\Vert \Exp_x^{-1}(y) \Vert$.
Parallel translation $\Gm_x^y$ denotes a the map which transports $v\in\mT_x\mM$ to $\Gm_x^y v \in\mT_y\mM$ along $\gm_{x\ra y}$ such that the vector stays constant by satisfying a zero-acceleration condition \citep[][equation (4.13)]{LeeJohn97}. 

For a smooth function $f:\mM\ra\mR$, $\grad f(x)\in\mT_x\mM$ denotes the Riemannian gradient of $f$ at $x\in\mM$ which satisfies $\frac{d}{dt}f(\gm(t)) = \langle\gm^\prime(t), \grad f(x)\rangle$~\citep[see][Sec 3.5.1 and (3.31)]{OptMan}. The Hessian of $f$ is defined jointly with the Riemannian structure of the manifold. The (directional) Hessian is
$
H(x)[\xi_x] := \nb_{\xi_x} \grad f, 
$
 and we use $H(x)[u,v]:=\langle u,H(x)[v]\rangle$ as a shorthand. We call $x\in \mM$ an $(\ep,-\sqrt{\rho\ep})$ saddle point when $\|\nb f(x)\|\le \ep$ and $\lb_{\min} (H(x))\le -\sqrt{\rho\ep}$.  We refer the interested reader to \citet{do2016differential} and \citet{LeeJohn97} which provide a thorough review on these important concepts of Riemannian geometry.
\vspace{-0.2cm}
\section{Perturbed Riemannian gradient algorithm}
\begin{figure}[htbp!]
      \centering
        \includegraphics[width = .4\textwidth]{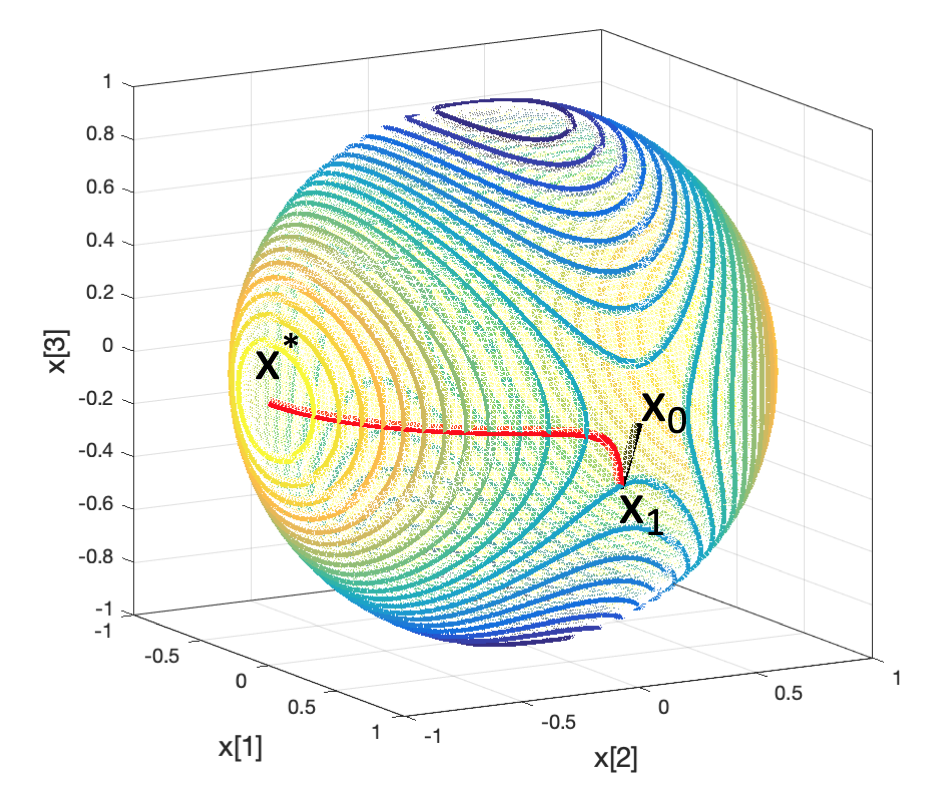}
      \caption{Function $f$ with saddle point on a sphere. $f(x) = x_1^2 - x_2^2 + 4x_3^2$. We plot the contour of this function on unit sphere. Algorithm~\ref{prgd} initializes at $x_0 = [1,0,0]$ (a saddle point), perturbs it towards $x_1$ and runs Riemannian gradient descent, and terminates at $x^* = [0,-1,0]$ (a local minimum). We amplify the first iteration to make saddle perturbation visible.}\label{fig_saddle}\vspace{-0.5cm}
    \end{figure}

\vspace{-0.2cm}
Our main Algorithm \ref{prgd} runs as follows: 
\vspace{-0.1cm}
\begin{enumerate}
\item Check the norm of the gradient: If it is large, do one step of Riemannian gradient descent, consequently the function value decreases. 
\item If the norm of gradient is small, it's either an approximate saddle point or a local minimum. Perturb the variable by adding an appropriate level of noise in its tangent space, map it back to the manifold and run a few iterations. 
\begin{enumerate}
\item If the function value decreases, iterates are escaping from the approximate saddle point (and the algorithm continues)
\item If the function value does not decrease, then it is an approximate local minimum (the algorithm terminates).
\end{enumerate}
\end{enumerate}
\vspace{-0.1cm}
\begin{algorithm}[htbp!]
\caption{Perturbed Riemannian gradient algorithm}\label{prgd}
\begin{algorithmic}
    \REQUIRE {Initial point $x_0 \in \mM$, parameters $\bt, \rho, K, \mathI$,
    accuracy $\ep$, probability of success $\delta$} (parameters defined in Assumptions \ref{assump:Lip_gra}, \ref{assump:Lip_Hes}, \ref{assump:sec_curv} and assumption of Theorem \ref{thm:main}). 
  \\
    Set constants:
    $\hat c\ge4$, $C:=C(K,\bt,\rho)$ (defined in Lemma \ref{lemma_2_pts} and proof of Lemma \ref{lemma2})\\
    ~~~~~and $\sqrt{c_{\max}}\le \frac{1}{56\hat c^2}$,
    $r=\frac{\sqrt{c_{\max}}}{\chi^2}\ep$,
    $\chi= 3\max\{ \log(\frac{d\bt(f(x_0) - f^*)}{\hat c\ep^2\dt}),4\}$.\\
    Set threshold values: $f_{\rm thres} = \frac{c_{\max}}{\chi^3}\sqrt{\frac{\ep^3}{\rho}}$, $g_{\rm thres} =  \frac{\sqrt{c_{\max}}}{\chi^2}\ep$, 
    $t_{\rm thres} = \frac{\chi}{c_{\max}}\frac{\bt}{\sqrt{\rho\ep}}$, $t_{\rm noise} = -t_{\rm thres}-1$.\\ 
    Set stepsize: $\eta=\frac{c_{\max}}{\bt}$.
    \WHILE {1}    
    \IF {$\|\grad f(x_t)\| \le g_{\rm thres}$ and $t-t_{\rm noise}>t_{\rm thres}$}
    \STATE {$t_{\rm noise}\la t$, $\tilde{x}_t \la x_t$, $x_t\la \Exp_{x_t}(\xi_t)$, $\xi_t$ uniformly sampled from $\mathbb{B}_{x_t}(r)\subset \mT_x\mM$.}
    \ENDIF
    \IF{$t-t_{\rm noise}=t_{\rm thres}$ and $f(x_t)-f(\tilde x_{t_{\rm noise}})>-f_{\rm thres}$}
     \STATE{ \textbf{output} {$\tilde x_{t_{\rm noise}}$} }
    \ENDIF
    \STATE {$x_{t+1}+ \la \Exp_{x_{t}} (- \min\{\eta, \frac{\mathI}{\|\grad f(x_{t})\|}\}\grad f(x_{t})).$}
    \STATE{$t\la t+1$.} 
   \ENDWHILE
     \end{algorithmic}
\end{algorithm}
Algorithm~\ref{prgd} relies on the manifold's exponential map, and is useful for cases where this map is easy to compute (true for many common manifolds). We refer readers to \citet[][pp. 81-86]{LeeJohn97} for the exponential map of sphere and hyperbolic manifolds, and \citet[][Example 5.4.2, 5.4.3]{OptMan} for the Stiefel and Grassmann manifolds. If the exponential map is not computable, the algorithm can use a 
retraction\footnote{A retraction is a first-order approximation of the exponential map which is often easier to compute.} instead, however our current analysis only covers the case of the exponential map. In Figure \ref{fig_saddle}, we illustrate a function with saddle point on sphere, and plot the trajectory of Algorithm~\ref{prgd} when it is initialized at a saddle point. 
    \vspace{-0.2cm}
\section{Main theorem: escape rate for perturbed Riemannian gradient descent}\label{main}
\vspace{-0.2cm}
We now turn to our main results, beginning with our assumptions and a statement of our main theorem. We then develop a brief proof sketch.

Our main result involves two conditions on function $f$ and one on the curvature of the manifold $\mM$.
\begin{assumption}[Lipschitz gradient] \label{assump:Lip_gra} There is a finite constant $\bt$ such that
\[
\|\grad f(y) - \Gm_x^y\grad f(x)\| \le \bt d(x,y) \quad \text{for all } x,y\in\mM.
\]
\end{assumption}
\begin{assumption}[Lipschitz Hessian] \label{assump:Lip_Hes} There is a finite constant $\rho$ such that
\[
\| H(y) - \Gm_x^y H(x) \Gm_y^x \|_2 \le \rho d(x,y) \quad \text{for all } x,y\in\mM.
\]
\end{assumption}
\begin{assumption}[Bounded sectional curvature] \label{assump:sec_curv} There is a finite constant $K$ such that 
\[
|K(x)[u,v]|\leq K \quad  \text{for all } x\in\mM \text{ and } u,v\in \mT_x\mM
\]
\end{assumption}
\vspace{-0.15cm}
 $K$ is an intrinsic parameter of the manifold capturing the curvature. We list a few examples here:
 (i) A sphere of radius $R$ has a constant sectional curvature $K = {1}/{R^2}$ \citep[][Theorem 1.9]{LeeJohn97}. If the radius is bigger, $K$ is smaller which means the sphere is less curved; (ii) A hyper-bolic space $H_R^n$ of radius $R$ has $K=-1/R^2$ \citep[][Theorem 1.9]{LeeJohn97}; (iii) For sectional curvature of the Stiefel and the Grasmann manifolds, we refer readers to \citet[][Section 5]{rapcsak2008sectional} and \citet{MR0229173}, respectively. 

Note that the constant $K$ is not directly related to the RLICQ parameter $R$ defined by \citet{Ge} which first requires describing the manifold by equality constraints. Different representations of the same manifold could lead to different curvature bounds, while sectional curvature is an intrinsic property of manifold. 
If the manifold is a sphere $\sum_{i=1}^{d+1} x_i^2 = R^2$, then $K = 1/R^2$, but more generally there is no simple connection. 
The smoothness parameters we assume are natural compared to some quantity from complicated compositions \citep[][Section 5.5]{lee2017first} or pullback \citep{zhang2018cubic,Criscitiello2019Efficiently}. With these assumptions, the main result of this paper is the following:
\begin{theorem}\label{thm:main}
Under Assumptions \ref{assump:Lip_gra},\ref{assump:Lip_Hes},\ref{assump:sec_curv}, let $C(K,\bt,\rho)$ be a function defined in Lemma \ref{lemma_2_pts}, $\hat \rho=\max\{\rho,C(K,\bt,\rho)\}$, if $\ep$ satisfies that
\begin{equation}\label{eq:eps_bd2}
        \ep\le\min\left\{ \frac{\hat\rho}{56\max\{c_2(K),c_3(K)\}\eta\bt}\log\left(\frac{d\bt}{\sqrt{\hat\rho\ep}\delta}\right), \left(\frac{\mathI\hat\rho}{12\hat c\sqrt{\eta\bt}}\log\left(\frac{d\bt}{\sqrt{\hat\rho\ep}\delta}\right)\right)^2\right\}
\end{equation}
where $c_2(K)$, $c_3(K)$ are defined in Lemma \ref{lem:Kar2}, 
then with probability $1-\delta$, perturbed Riemannian gradient descent with step size $c_{\max}/\bt$ converges to a $(\ep, -\sqrt{\hat \rho\ep})$-stationary point of $f$ in 
\begin{equation*}
O\Bigg(\frac{\bt(f(x_0) - f(x^*))}{\ep^2}\log^4\bigg(\frac{\bt d(f(x_0) - f(x^*))}{\ep^2\delta}\bigg)\Bigg)
\end{equation*}
iterations. 
\end{theorem}
{\bf Proof roadmap.}
For a function satisfying smoothness condition (Assumption \ref{assump:Lip_gra} and \ref{assump:Lip_Hes}), we use a local upper bound of the objective based on the third-order Taylor expansion (see supplementary material Section A for a review),
\begin{equation*}
    \begin{split}
        f(u) &\le f(x) + \langle\grad f(x), \Exp_x^{-1}(u)\rangle + \frac{1}{2} H(x)[\Exp_x^{-1}(u),\Exp_x^{-1}(u)] + \frac{\rho}{6}\|\Exp_x^{-1} (u)\|^3.
    \end{split}
\end{equation*}
When the norm of the gradient is large (not near a saddle), the following lemma guarantees the decrease of the objective function in one iteration.
\begin{lemma}\label{lem:first_order_progress}\citep{manifold_second}
Under Assumption~\ref{assump:Lip_gra}, by choosing $\bar\eta=\min\{\eta, \frac{\mathI}{\|\grad f(u)\|}\} = O(1/\bt)$, the Riemannian gradient descent algorithm is monotonically descending, $f(u^+) \!\le \!f(u) \!-\! \frac{1}{2}\bar\eta \|\grad f(u)\|^2$.
\end{lemma}


Thus our main challenge in proving the main theorem is the Riemannian gradient behaviour at an approximate saddle point:

 1. Similar to the Euclidean case studied by \citet{GeSham}, we need to bound the ``thickness'' of the ``stuck region'' where the perturbation fails. We still use a pair of hypothetical auxiliary sequences and study the ``coupling'' sequences. When two perturbations couple in the thinnest direction of the stuck region, their distance grows and one of them escapes from saddle point. 
 
2. However our iterates are evolving on a manifold rather than a Euclidean space, so our strategy is to map the iterates back to an appropriate fixed tangent space where we can use the Euclidean analysis. This is done using the inverse of the exponential map and various parallel transports.

3. Several key challenges arise in doing this. Unlike \citet{GeSham}, the structure of the manifold interacts with the local approximation of the objective function in a complicated way. On the other hand, unlike recent work on Riemannian optimization by \citet{Rie_trust_region}, we do not have access to a second order oracle and we need to understand how the sectional curvature and the injectivity radius (which both capture intrinsic manifold properties) affect the behavior of the first order iterates.

4. Our main contribution is to carefully investigate how the various approximation errors arising from (a) the linearization of the iteration couplings and (b) their mappings to a common tangent space can be handled on manifolds with bounded sectional curvature. We address these challenges in a sequence of lemmas (Lemmas \ref{lem:Kar1} through \ref{lem:par}) we combine to   linearize the coupling iterations in a common tangent space and precisely control the approximation error. 
This result is formally stated in the following lemma.

\begin{lemma}\label{lemma_2_pts} 
Define $\gm = \sqrt{\hat\rho\ep}$, $\kappa = \frac{\bt}{\gm}$, and $\mathS = \sqrt{\eta\bt} \frac{\gm}{\hat\rho} \log^{-1}(\frac{d\ka}{\delta})$. 
Let us consider $x$ be a $(\ep,-\sqrt{\hat\rho\ep})$ saddle point, and  define $u^+ = \Exp_u(-\eta \grad f(u))$ and $w^+ = \Exp_w(-\eta \grad f(w))$.  Under Assumptions~\ref{assump:Lip_gra},~\ref{assump:Lip_Hes},~\ref{assump:sec_curv}, 
if all pairwise distances between $u,w,u^+,w^+,x$ are less than $12\mathS$, then for some explicit constant  $C(K,\rho,\bt)$ depending only on $K,\rho,\bt$, there is
\begin{align}
 &\quad \|\Expn_x(w^+) - \Expn_x(u^+) -(I-\eta H(x))(\Expn_x(w) - \Expn_x(u))\| \label{eq:exp_inv_bound}\\
 &\le C(K,\rho,\bt) d(u,w)\lc d(u,w)+d(u,x) + d(w,x)\rc. \notag
\end{align}
 \end{lemma}

 The proof of this lemma includes novel contributions by strengthen known result (Lemmas \ref{lem:Kar1}) and also combining known inequalities in novel ways (Lemmas \ref{lem:Kar2}~to~\ref{lem:par}) that allow us to control all the approximation errors and arrive at the tight rate of escape for the algorithm.

\vspace{-0.2cm}
\section{Proof of Lemma \ref{lemma_2_pts}}\label{uvw}
\vspace{-0.2cm}
Lemma \ref{lemma_2_pts} controls the error of the linear approximation of the iterates when mapped in $T_x\mM$. 
In this section, we assume that all points are within a region of diameter $R:= 12\mathS\le\mathI$ (inequality follows from \eq{eps_bd2}
), i.e., the distance of any two points in the following lemmas are less than $R$.
The proof of Lemma~\ref{lemma_2_pts} is based on the sequence of following lemmas.
\begin{lemma}\label{lem:Kar1} 
 Let $x\in \mM$ and $y,a\in T_x\mM$. Let us denote by $z=\Exp_x(a)$ then under Assumption~\ref{assump:sec_curv}
\begin{equation}
\begin{split}
d(\Exp_x(y+a), \Exp_{z}(\Gamma_x^{z}y))\leq c_1(K) \min\{\Vert a \Vert, \Vert y \Vert\} (\Vert a \Vert+ \Vert y \Vert)^2. 
\end{split}
\end{equation}
\end{lemma}
This lemma tightens the result of \citet[][C2.3]{karcher}, which only shows an upper-bound 
$O(\Vert a \Vert(\Vert a \Vert+ \Vert y \Vert)^2)$.
We prove the upper-bound 
$O(\Vert y \Vert(\Vert a \Vert+ \Vert y \Vert)^2)$  in the supplement. 
\begin{figure}[t!]
      \centering
      \subfloat[]{
        \includegraphics[width = .3\textwidth]{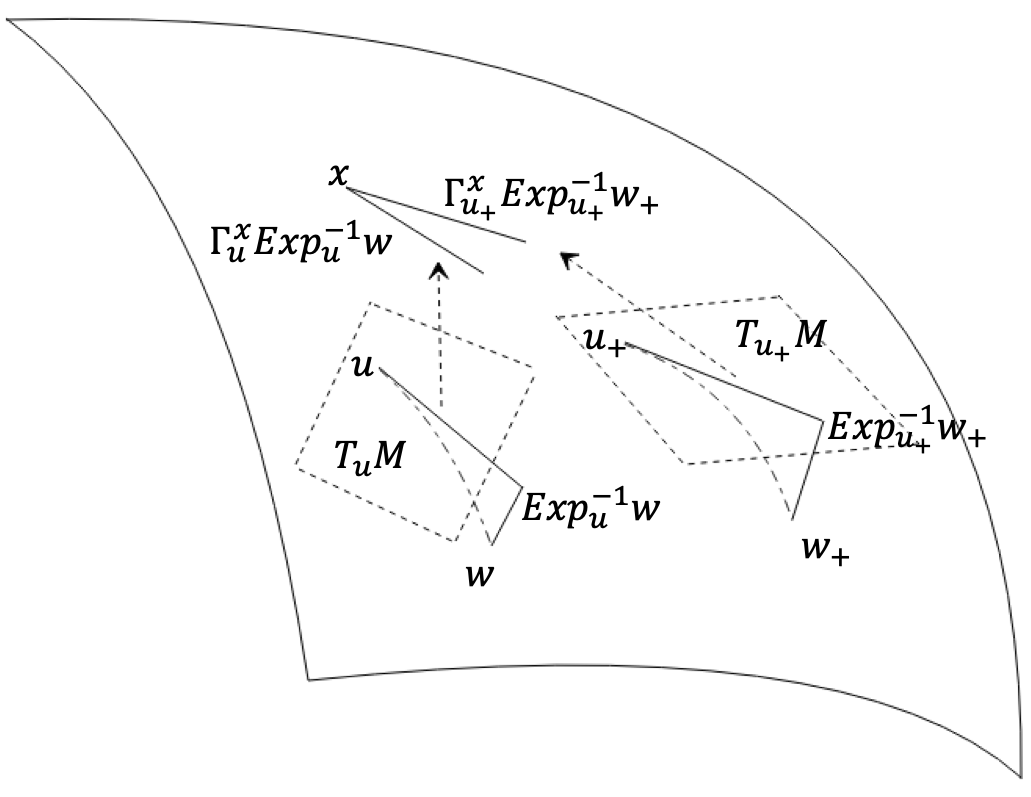}\label{fig1}
      }
      \subfloat[]{
        \includegraphics[width = .3\textwidth]{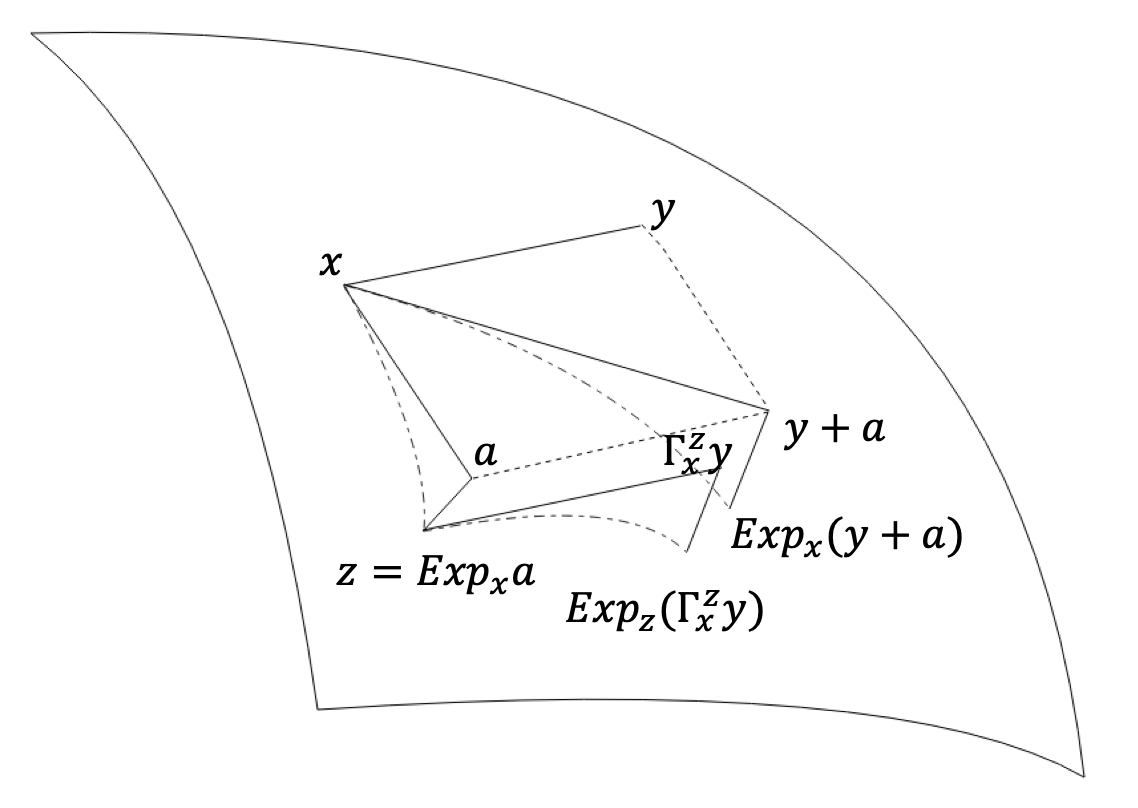}\label{fig2}
      }
      \caption{(a) \eq{fig_1}. First map $w$ and $w_+$ to $\mT_u\mM$ and $\mT_{u_+}\mM$, and transport the two vectors to $\mT_x\mM$, and get their relation. (b) Lemma \ref{lem:Kar1} bounds the difference of two steps starting from $x$: (1) take $y+a$ step in $\mT_x\mM$ and map it to manifold, and (2) take $a$ step in $\mT_x\mM$, map to manifold, call it $z$, and take $\Gm_x^z y$ step in $\mT_x\mM$, and map to manifold. $\Exp_z(\Gm_x^z y)$ is close to $\Exp_x(y+a)$.}\vspace{-0.2cm}
    \end{figure}
We also need the following lemma showing that both the exponential map  and its inverse are Lipschitz.
\begin{lemma}\label{lem:Kar2} Let $x,y,z\in M$, and the distance of each two points is no bigger than $R$. Then under assumption~\ref{assump:sec_curv}
\begin{equation*}
    (1+c_2(K)R^2)^{-1}d(y,z)\leq \Vert \Exp_x^{-1}(y)-\Exp_x^{-1}(z)\Vert \leq(1+c_3(K)R^2)d(y,z).
\end{equation*}
\end{lemma}  
Intuitively this lemma relates the norm of the difference of two vectors of $\mT_x\mM$ to the distance between the corresponding points on the manifold $\mM$ and follows from bounds on the Hessian of the square-distance function~\citep[][Ex. 4 p. 154]{Sak96}.
The upper-bound is directly proven by \citet[][Proof of Cor. 1.6]{karcher}, and we prove the lower-bound via Lemma~\ref{lem:Kar1} in the supplement.

The following contraction result is fairly classical and is proven using the Rauch comparison theorem from differential geometry~\citep{MR2394158}.
\begin{lemma}\label{lem:dis}
\citep[][Lemma 1]{rapidmix} Under Assumption~\ref{assump:sec_curv}, for $x,y\in \mM$ and $w\in T_x\mM$,
\[
d(\Exp_x(w),\Exp_y(\Gamma_x^y w))\leq c_4(K) d(x,y).
\]
\end{lemma}
Finally we need the following corollary of the Ambrose-Singer theorem~\citep{ambrosesinger}.
\begin{lemma}\label{lem:par}
\citep[][Section 6]{karcher} Under Assumption~\ref{assump:sec_curv}, for $x,y,z\in \mM$ and $w\in T_x\mM$,
\[
\Vert \Gamma_y^{z}\Gamma_{x}^yw -\Gamma_{x}^{z}w \Vert \leq c_5(K) d(x,y)d(y,z) \Vert w\Vert.
\]
\end{lemma}
Lemma \ref{lem:Kar1} through \ref{lem:par} are mainly proven in the literature, and we make up the missing part in Supplementary material Section B. Then we prove Lemma \ref{lemma_2_pts} in Supplementary material Section B.

The spirit of the proof is to linearize the manifold using the exponential map and its inverse, and to carefully bounds the various error terms caused by the approximation. Let us denote by \mbox{$\theta= d(u,w)+d(u,x) + d(w,x)$}. 

1. We first show using twice Lemma~\ref{lem:Kar1} and Lemma~\ref{lem:dis} that 
\[
d(\Exp_u(\Exp_u^{-1}(w)-\eta \Gamma_w^u \grad f(w)), \Exp_u(-\eta \grad f(u)+\Gamma_{u_+}^u \Exp_{u_+}^{-1}(w_+)) ) =O(\theta  d(u,w)).\]
2. We use Lemma~\ref{lem:Kar2} to linearize this iteration in $\mT_u\mM$ as
\[
 \Vert \Gamma_{u_+}^u \Exp_{u_+}^{-1}(w_+) -\Exp_u^{-1}(w) +\eta [\grad f(u)-\Gamma_w^u \grad f(w)] \Vert = O(\theta d(u,w)). 
\]
3. Using the Hessian Lipschitzness  
\[
 \Vert \Gamma_{u_+}^u \Exp_{u_+}^{-1}(w_+)) -\Exp_u^{-1}(w) +\eta H(u) \Exp_u^{-1}(w) \Vert = O(\theta d(u,w)).
\]
3. We use Lemma~\ref{lem:par}  to map to $T_{x}\mM$ and the Hessian Lipschitzness to compare $H(u)$ to $H(x)$. This is an important intermediate result (see Lemma~1 in Supplementary material Section~B).
\begin{align}
   \Vert \Gamma_{u_+}^x \Exp_{u_+}^{-1}(w_+) -\Gamma_{u}^x\Exp_u^{-1}(w) +\eta H(x)\Gm_u^x \Exp^{-1}_u(w) \Vert = O(\theta d(u,w)).  \label{eq:fig_1}
\end{align}

4. We use Lemma \ref{lem:Kar1} and \ref{lem:Kar2} to approximate two iteration updates in $\mT_x\mM$. 
\begin{align}\label{eq:compaa}
   \|\Expn_x(w)- (\Expn_x(u) + \Gm_u^x\Expn_u(w)) \|\le O(\theta d(u,w)). 
\end{align}
And same for the $u_+,w_+$ pair replacing $u,w$.

5. Combining \eq{fig_1} and \eq{compaa} together, we obtain
\[
\|\Expn_x(w^+) - \Expn_x(u^+) -(I-\eta H(x))(\Expn_x(w) - \Expn_x(u))\| \le O(\theta d(u,w)). 
\]
Now note that, the iterations  $u,u_+,w,w_+$ of the algorithm are both on the manifold. We use $\Expn_x(\cdot)$ to map them to the same tangent space at $x$. 

Therefore we have linearized the two coupled trajectories $\Expn_x(u_t)$ and $\Expn_x(w_t)$ in a common tangent space, and  we can modify the Euclidean escaping saddle analysis thanks to the error bound we proved in Lemma~\ref{lemma_2_pts}. 

\vspace{-0.2cm}
\section{Proof of main theorem}\label{comp_pr}
\vspace{-0.2cm}
In this section we suppose all assumptions in Section~\ref{main} hold. The proof strategy is to show with high probability that the function value decreases of $\mF$ in $\mathT$ iterations at an approximate saddle point. Lemma \ref{lemma1} suggests that, if after a perturbation and $\mathT$ steps, the iterate is $\Omega(\mathS)$ far from the approximate saddle point, then the function value decreases. If the iterates do not move far, the perturbation falls in a stuck region. Lemma \ref{lemma2} uses a coupling strategy, and suggests that the width of the stuck region is small in the negative eigenvector direction of the Riemannian Hessian. 

Define
\begin{align*}
\mF = \eta\bt \frac{\gm^3}{\hat\rho^2} \log^{-3}(\frac{d\ka}{\delta}), \ \mG = \sqrt{\eta\bt} \frac{\gm^2}{\hat\rho} \log^{-2}(\frac{d\ka}{\delta}),\
\mathT = \frac{\log(\frac{d\ka}{\delta})}{\eta\gm}.
\end{align*}
At an approximate saddle point $\tilde x$, 
let $y$ be in the neighborhood of $\tilde x$ where $d(y,\tilde x)\le \mathI$, denote
\begin{equation*}
\tilde f_y(x) := f(y) + \langle \grad f(y), \Exp^{-1}_{y}(\tilde x)\rangle + \frac{1}{2} H(\tilde x)[\Exp^{-1}_{y}(\tilde x), \Exp^{-1}_{y}(\tilde x)].
\end{equation*}
Let $\| \grad f(\tilde x)\|\le \mG$ and $\lb_{\min} (H(\tilde x)) \le -\gm$. We consider two iterate sequences, $u_0,u_1,...$ and $w_0,w_1,...$ where $u_0,w_0$ are two perturbations at $\tilde x$. 
\begin{lemma}\label{lemma1}
Assume Assumptions \ref{assump:Lip_gra}, \ref{assump:Lip_Hes}, \ref{assump:sec_curv} and \eq{eps_bd2} hold. There exists a constant $c_{\max}$, $\forall \hat c>3, \dt\in (0,\frac{d\ka}{e}]$, for any $u_{0}$ with $d(\tilde x, u_0)\le2\mathS/(\ka \log(\frac{d\ka}{\dt}))$, $\ka = \bt/\gm$.
\begin{equation*}
  T = \min\lbig \inf_t\lbig t|\tilde f_{u_0}(u_t) - f(u_0) \le -3\mF\rbig, \hat c \mathT  \rbig,
  \end{equation*}
  then $\forall \eta\le c_{\max}/\bt$, we have $\forall 0<t<T$, $d(u_0, u_t)\le 3(\hat c\mathS)$.
\end{lemma}
\begin{lemma}\label{lemma2}
Assume Assumptions \ref{assump:Lip_gra}, \ref{assump:Lip_Hes}, \ref{assump:sec_curv} and \eq{eps_bd2} hold. Take two points $u_0$ and $w_0$ which are perturbed from an approximate saddle point, where $d(\tilde x, u_0)\le 2\mathS/(\ka \log(\frac{d\ka}{\dt}))$, 
$\Expn_{\tilde x}(w_0) - \Expn_{\tilde x}(u_0) = \mu r e_1$,
$e_1$ is the smallest eigenvector\footnote{``smallest eigenvector'' means the eigenvector corresponding to the smallest eigenvalue.} of $H(\tilde x)$, $\mu\in[\dt/(2\sqrt{d}),1]$, and the algorithm runs two sequences $\{u_t\}$ and $\{w_t\}$ starting from $u_0$ and $w_0$. Denote
\begin{equation*}
  T = \min\lbig \inf_t\lbig t|\tilde f_{w_0}(w_t) - f(w_0) \le -3\mF\rbig, \hat c \mathT  \rbig,
\end{equation*}
  then $\forall \eta\le c_{\max}/l$, if $\forall 0<t<T$, $d(\tilde x, u_t)\le 3(\hat c\mathS)$, we have $T<\hat c \mathT$.
\end{lemma}
We prove Lemma \ref{lemma1} and \ref{lemma2} in supplementary material Section C.  We also prove, in the same section, the main theorem using the  coupling strategy of \citet{GeSham}. 
but with the additional difficulty of taking into consideration the effect of the Riemannian geometry (Lemma \ref{lemma_2_pts}) and the injectivity radius. 


\section{Examples}
\begin{figure}[tbp!]
      \centering
      \subfloat[]{
       \includegraphics[width = 4.5cm,height=3cm]{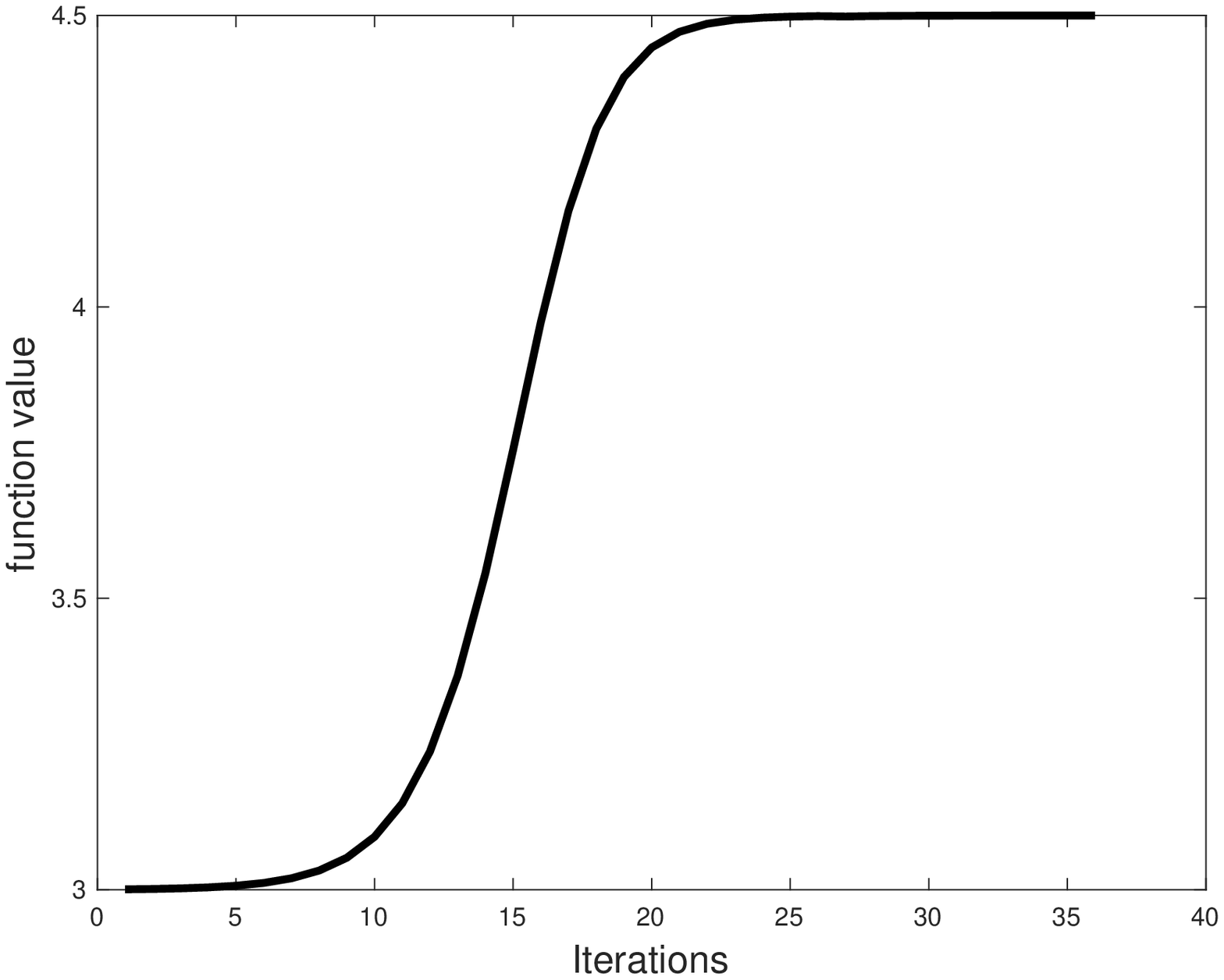}\label{f1}
      }
      \qquad
      \subfloat[]{
        \includegraphics[width = 4.5cm,height=3cm]{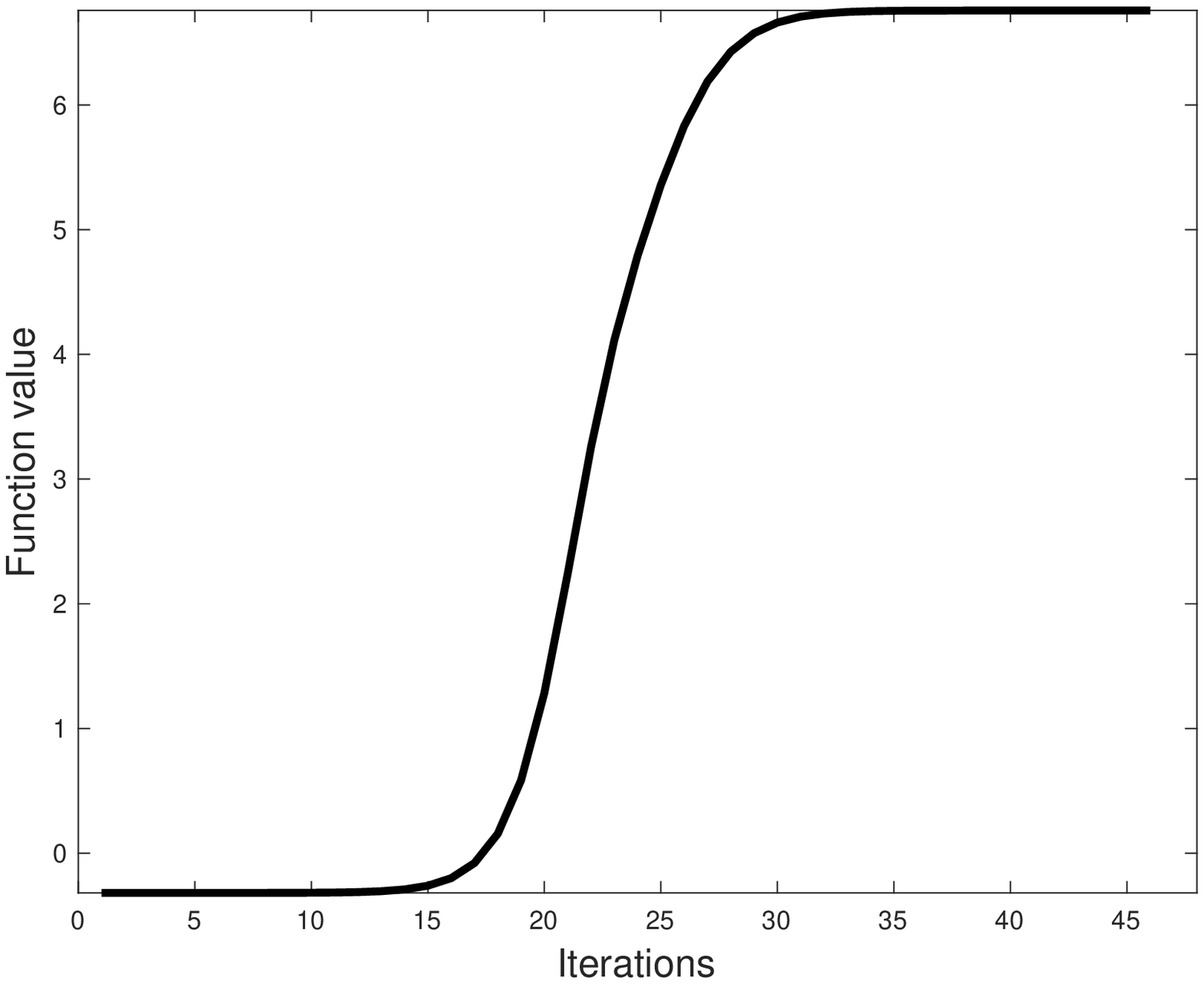}\label{f2}
      }
      \caption{(a) kPCA problem. 
      We start from an approximate saddle point, and it converges to a local minimum (which is also global minimum). (b) Burer-Monteiro approach 
      Plot $f(Y) = \half \mathrm{trace}(AYY^T)$ versus iterations. We start from the saddle point, and it converges to a local minimum (which is also global minimum).}\label{fig:exp}\vspace{-0.2cm}
\end{figure}

\vspace{-0.2cm}
\paragraph{kPCA.}
We consider the kPCA problem, where we want to find the $k\leq n$ principal eigenvectors of a symmetric matrix $H\in\mathbb{R}^{n\times n}$,  as an example~\citep{manifold_first}. This corresponds to  
\begin{align*}
 \min_{X\in\mathbb{R}^{n\times k}} \ -\half \mathrm{tr}(X^THX)\quad 
 \mbox{subject to} \ X^TX = I, 
\end{align*}
which is an optimization problem on the Grassmann manifold defined by the constraint $X^TX = I$. If the eigenvalues of $H$ are distinct, we denote by $v_1$,...,$v_n$ the eigenvectors of $H$, corresponding to eigenvalues with decreasing order. Let
$
    V^* = [v_1,...,v_k]
$
be the matrix with columns composed of the top $k$ eigenvectors of $H$, then the local minimizers of the objective function are $V^*G$ for all unitary matrices $G\in\mR^{k\times k}$. Denote also by
$
    V = [v_{i_1},...,v_{i_k}]
$
the matrix with columns composed of $k$ distinct eigenvectors, then the first order stationary points of the objective function (with Riemannian gradient being $0$) are $VG$ for all unitary matrices $G\in\mR^{k\times k}$. In our numerical experiment, we choose $H$ to be a diagonal matrix $H=\diag(0,1,2,3,4)$ and let $k=3$. The Euclidean basis $(e_i)$ are an eigenbasis of $H$ and the first order stationary points of the objective function are $[e_{i_1},e_{i_2},e_{i_3}]G$ with distinct basis and $G$ being unitary. The local minimizers are $[e_3,e_4,e_5]G$. We start the iteration at $X_0 = [e_2,e_3,e_4]$ and see in Fig.~\ref{fig:exp} the algorithm converges to a local minimum. 
\vspace{-0.1cm}
\paragraph{Burer-Monteiro approach for certain low rank problems.}
Following \citet{BM}, we consider, for $A\in\mathbb{S}^{d\times d}$ and $r(r+1)/2\le d$, the problem 
\begin{align*}
\minl{X\in\mathbb{S}^{d\times d}} \trace(AX), \ s.t.\ \diag(X)=1,X\succeq0, \rank(X)\le r.
\end{align*}
We factorize $X$ by $YY^T$ with an overparametrized $Y\in\mathbb{R}^{d\times p}$ and $p(p+1)/2\ge d$. Then any local minimum of
\begin{align*}
\minl{Y\in\mathbb{R}^{d\times p}} \trace(AYY^T), \ s.t.\ \diag(YY^T)=1,
\end{align*}
is a global minimum where $YY^T=X^*$ \citep{BM}. Let $f(Y) = \frac{1}{2}\trace (AYY^T)$. In the experiment, we take $A\in\mathbb{R}^{100\times 20}$ being a sparse matrix that only the upper left $5\times 5$ block is random and other entries are $0$. Let the initial point $Y_0\in\mathbb{R}^{100\times 20}$, such that $  (Y_0)_{i,j}=1 $ for $ 5j-4\le i \le 5j$ and $  (Y_0)_{i,j}=0 $ otherwise.
Then $Y_0$ is a saddle point. 
We see in Fig.~\ref{fig:exp} the algorithm converges to the global optimum.

\paragraph{Summary}
We have shown that for the constrained optimization problem of minimizing $f(x)$ subject to a manifold constraint as long as the function and the manifold are appropriately smooth, a perturbed Riemannian gradient descent algorithm will escape saddle points with a rate of order $1/\epsilon^2$ in the accuracy $\epsilon$, polylog in manifold dimension $d$, and depends polynomially on the curvature and smoothness parameters.

A natural extension of our result is to consider other variants of gradient descent, such as the heavy ball method, Nesterov's acceleration, and the stochastic setting. The question is whether these algorithms with appropriate modification (with manifold constraints) would have a fast convergence to second-order stationary point (not just first-order stationary as studied in recent literature), and whether it is possible to show the relationship between convergence rate and smoothness of manifold.
\bibliography{bibfile}
\bibliographystyle{icml2019}
\appendix

\section*{Appendix} 

\subsection*{Organization of the Appendix}

In Appendix~\ref{app:taylor} we review classical results on the Taylor expansions for functions on Riemannian manifold. In Appendix~\ref{app:mani} we provide the proof of Lemma~\ref{lemma_exp_inv_bound} which requires to expand the iterates on the tangent space in the the saddle point. Finally, in Appendix~\ref{app:proof}, we provide the proofs of Lemma~\ref{lemma1_a} and Lemma~\ref{lemma2_a} which enable to prove the main theorem of the paper.

Throughout the paper we assume that the objective function and the manifold are smooth. Here we list the assumptions that are used in the following lemmas.
\setcounter{assumption}{0}
\begin{assumption}[Lipschitz gradient] \label{assump:Lip_gra_a} There is a finite constant $\bt$ such that
\[
\|\grad f(y) - \Gm_x^y\grad f(x)\| \le \bt d(x,y) \quad \text{for all } x,y\in\mM.
\]
\end{assumption}
\begin{assumption}[Lipschitz Hessian] \label{assump:Lip_Hes_a} There is a finite constant $\rho$ such that
\[
\| H(y) - \Gm_x^y H(x) \Gm_y^x \|_2 \le \rho d(x,y) \quad \text{for all } x,y\in\mM.
\]
\end{assumption}
\begin{assumption}[Bounded sectional curvature] \label{assump:sec_curv_a} There is a finite constant $K$ such that 
\[
|K(x)[u,v]|\leq K \quad  \text{for all } x\in\mM \text{ and } u,v\in \mT_x\mM
\]
\end{assumption}

\section{Taylor expansions on Riemannian manifold}\label{app:taylor}

We provide here the Taylor expansion for functions and gradients of functions defined on a Riemannian manifold.

\subsection{Taylor expansion for the gradient}\label{Tg}
For any point $x\in\mM$ and $z\in\mM$ be a point in the neighborhood of $x$ where the geodesic $\gm_{x\ra z}$ is defined. 
\begin{align}
 \Gm_z^x( \grad f(z))
&= \grad f(x) + \nb_{\gm^\prime_{x\ra z}(0)} \grad f
\textstyle +\int_0^1 ( \Gm_{\gm_{x\ra z}(\tau)}^x\nb_{ \gm^\prime_{x\ra z}(\tau)} \grad f - \nb_{\gm^\prime_{x\ra z}(0)} \grad f) d\taux\notag\\
&=  \grad f(x) + \nb_{\gm^\prime_{x\ra z}(0)} \grad f + \Dt(z), \label{eq:Taylor}
\end{align}
where $ \Dt(z):= \textstyle\int_0^1 {( \Gm_{\gm_{x\ra z}(\tau)}^x \nb_{\gm^\prime_{x\ra z}(\tau)} \grad f} 
- \nb_{\gm^\prime_{x\ra z}(0)} \grad f) d\tau$.
The Taylor approximation in \eq{Taylor} is proven by \citet[][Lemma 7.4.7]{OptMan}.

\subsection{Taylor expansion for the function}\label{Tf}
Taylor expansion of the gradient enables us to approximate the iterations of the main algorithm, but obtaining the convergence rate of the algorithm requires proving that the function value decreases following the iterations. We need to give the Taylor expansion of $f$ with the parallel translated gradient on LHS of  \eq{Taylor}. To simplify the notation, let $\gm$ denote the $\gm_{x\ra z}$.
\begin{subequations}
\begin{align}
f(z) \!-\! f(x)\! &= \int_0^1{ \frac{d}{d\tau} f(\gm(\tau)) d\tau}\\
	&=  \int_0^1{ \langle\gm^\prime(\tau), \grad f(\gm(\tau))\rangle d\tau}\\
	&=\!\int_0^1\!\!\!\! {\langle\Gm_{\gm(\tau)}^x\!\gm^\prime(\!\tau\!), \Gm_{\gm(\tau)}^x\grad f(\gm(\tau))\rangle d\tau}\label{eq:trans_preserve_ip}\\
	&= \int_0^1 {\langle\gm^\prime(0), \Gm_{\gm(\tau)}^0\grad f(\gm(\tau))\rangle d\tau}\label{eq:trans_speed}\\
	&= \int_0^1 {\langle\gm^\prime(0), \grad f(x) + \nb_{\tau \gm^\prime(0)} \grad f + \Dt(\gm(\tau)) \rangle d\tau}\label{eq:Taylorf_mani1}\\
	&=\textstyle \langle\gm^\prime(0), \grad f(x) + \frac{1}{2}\nb_{\gm^\prime(0)} \grad f +\bar\Dt(z) \rangle. 
\end{align}\label{eq:Taylorf_mani}
\end{subequations} 
$\Dt(z)$ is defined in \eq{Taylor}. $\bar\Dt(z) = \int_0^1 \Dt(\gm(\tau)) d\tau$. The second line is just rewriting by definition. \eq{trans_preserve_ip} means the parallel translation preserves the inner product \citep[][Prop. 14.16]{DiffGeo}. \eq{trans_speed} uses $\Gm_{\gm(t)}^x\gm^\prime(t) = \gm^\prime(0)$, meaning that the velocity stays constant along a geodesic \citep[][(5.23)]{OptMan}. \eq{Taylorf_mani1} uses \eq{Taylor}. In Euclidean space, the Taylor expansion is 
\begin{align}
 f(z) - f(x)= \langle z, \nb f(x) + \nb^2f(x)z + \int_0^1{(\nb^2f(\tau z) -\nb^2f(x))zd\tau} \rangle. \label{eq:Taylorf}
\end{align} 
Compare \eq{Taylorf_mani} and \eq{Taylorf}, $z$ is replaced by $\gm^\prime(0) := \gm_{x\ra z}^\prime(0)$ and $\tau z$ is replaced by $\tau\gm_{x\ra z}^\prime(0)$ or $\gm_{x\ra z}(\tau)$.

Now we have
\begin{equation*}
f(u_t) = f(x) + \langle \gm^\prime(0), \grad f(x) \rangle +\frac{1}{2}H(x)[\gm^\prime(0), \gm^\prime(0)] +   \langle \gm^\prime(0),\bar \Dt(u_t) \rangle. 
\end{equation*}

\section{Linearization of the iterates in a fixed tangent space  }\label{app:mani}

In this section we linearize the progress of the iterates of our algorithm in a fixed tangent space $\mT_x\mM$. We always assume here that all points are within a region of diameter $R:= 12\mathS\le\mathI$. In the course of the proof we need several auxilliary lemmas which are stated in the last two subsections of this section.

\subsection{Evolution of $\Exp_u^{-1}(w)$}
\begin{figure}[t!]
      \centering
        \includegraphics[width = .5\textwidth]{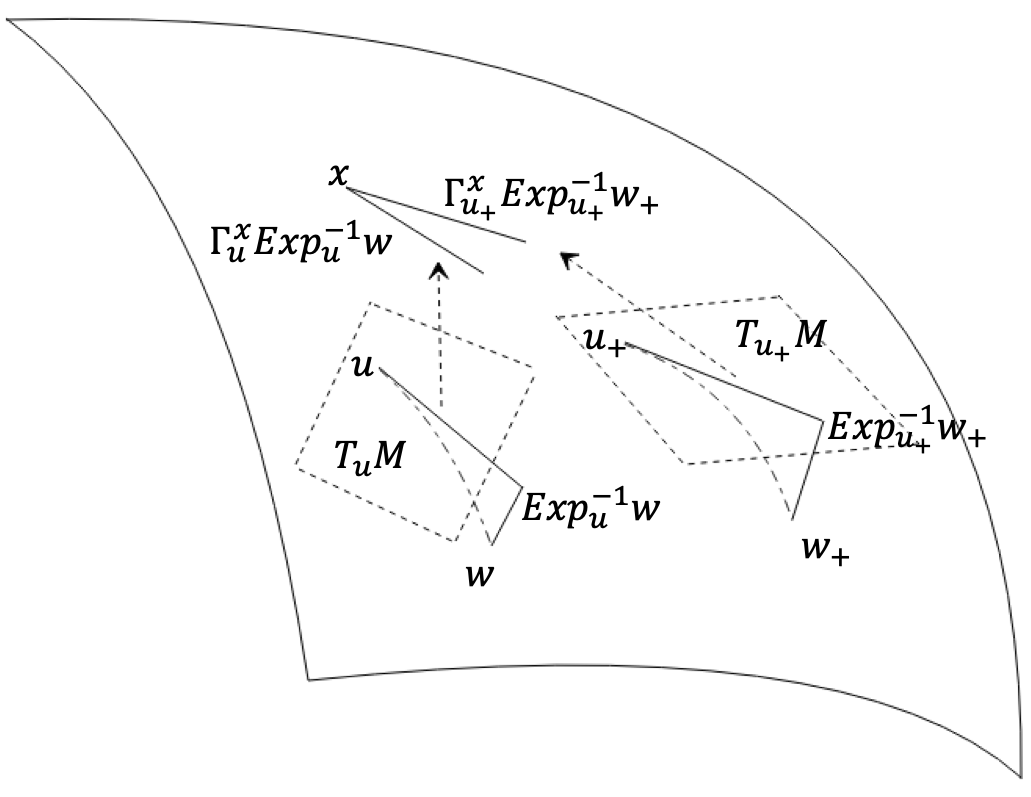}
      \caption{Lemma~\ref{lemma_2_pts_a}. First map $w$ and $w_+$ to $\mT_u\mM$ and $\mT_{u_+}\mM$, and transport the two vectors to $\mT_x\mM$, and get their relation.}\label{fig1}\vspace{-0.2cm}
    \end{figure}
We first consider the evolution of $\Exp_u^{-1}(w)$ in  a fixed tangent space $\mT_x\mM$. We show in the following lemma that it approximately follows a linear reccursion. 
 \begin{lemma}\label{lemma_2_pts_a}
Define $\gm = \sqrt{\hat\rho\ep}$, $\kappa = \frac{\bt}{\gm}$, and $\mathS = \sqrt{\eta\bt} \frac{\gm}{\hat\rho} \log^{-1}(\frac{d\ka}{\delta})$. 
Let us consider $x$ be a $(\ep,-\sqrt{\hat\rho\ep})$ saddle point, and  define $u^+ = \Exp_u(-\eta \grad f(u))$ and $w^+ = \Exp_w(-\eta \grad f(w))$.  Under Assumptions~\ref{assump:Lip_gra_a},~\ref{assump:Lip_Hes_a},~\ref{assump:sec_curv_a}, 
if all pairwise distances between $u,w,u^+,w^+,x$ are less than $12\mathS$, then for some explicit constant  $C_1(K,\rho,\bt)$ depending only on $K,\rho,\bt$, there is
 \begin{equation*}
 \begin{split}
      & \quad \|\Gm_{u^+}^x \Exp^{-1}_{u^+}(w^+) - (I - \eta H(x))\Gm_u^x \Exp^{-1}_u(w) \| \\
     &\le C_1(K,\rho,\bt) d(u,w)\lc d(u,w)+d(u,x) + d(w,x)\rc.
 \end{split}
 \end{equation*}
 for some explicit function $C_1$.
\end{lemma}
This lemma is illustrated in Fig.~\ref{fig1}.

\begin{proof}
Denote $-\eta \grad f(u) = v_u$, $-\eta \grad f(w) = v_w$. $v$ is a smooth map. We first prove the following claim.
\begin{claim}
\[
d(u_+,w_+)\leq  c_6(K)d(u,w),
\]
where $c_6(K)=c_4(K)+1+c_2(K)R^2$.
\end{claim}
To show this, note that
\[
d(u_+,w_+)\leq d(u_+,\tilde w_+)+d(\tilde w_+,w_+),
\]
and using Lemma \ref{lem:dis_a} with $\tilde{w}_+ = \Exp_w(\Gamma_u^w v_u)$,
\begin{equation*}
\begin{split}
d(\tilde w_+,w_+)&= d( \Exp_w(v_w), \Exp_w(\Gamma_u^w v_u))\\
&\leq (1+c_2(K)R^2)\Vert v_w-\Gamma_u^w v_u\Vert \\
&\leq \bt (1+c_2(K)R^2)d(u,w).
\end{split}
\end{equation*}
Using Lemma \ref{lem:dis_a}, 
\begin{equation}
\begin{split}
d(\tilde w_+,u_+)\le c_4(K)d(u,w).
\end{split}
\end{equation}
Adding the two inequalities proves the claim. 

We use now Lemma \ref{lem:Kar1_a} between $(u,w,u_+,w_+)$ in two different ways. First let us use it for $a=\Exp_u^{-1}(w)$ and $y=\Gamma_w^u v_w$. We obtain:
\begin{equation}
d(w_+,\Exp_u(\Exp_u^{-1}(w)+\Gamma_w^u v_w))\leq c_1(K) d(u,w) (d(u,w)^2+\Vert v_w\Vert^2). 
\end{equation}
Then we use it for $a=\Exp_u^{-1}(v_u)$ and $y=\Gamma_{u_+}^u \Exp_{u_+}^{-1}(w_+)$ which yields
\begin{align*}
 &\quad d(w_+,\Exp_u(v_u+\Gamma_{u_+}^u \Exp_{u_+}^{-1}(w_+))) \\
 &\leq c_1(K) d(u_+,w_+) ( d(u_+,w_+)^2+\Vert v_u\Vert^2)\\
 &\leq c_1(K) c_5(K,\Vert v_u\Vert, \Vert v_w\Vert)d(u,w)\cdot\Big[ c_5(K,\Vert v_u\Vert, \Vert v_w\Vert)^2d(u,w)^2+\Vert v_u\Vert^2\Big].
\end{align*}
Using the triangular inequality we have
\begin{align*}
 &\quad d(\Exp_u(\Exp_u^{-1}(w)+\Gamma_w^u v_w), \Exp_u(v_u+\Gamma_{u_+}^u \Exp_{u_+}^{-1}(w_+)) ) \\
 &\leq d(w_+,\Exp_u(\Exp_u^{-1}(w)+\Gamma_w^u v_w))   + d(w_+,\Exp_u(v_u+\Gamma_{u_+}^u \Exp_{u_+}^{-1}(w_+))) \\
 &\leq  c_7 d(u,w)
 \end{align*}
 with $c_7$ defined as
 \begin{align*}
    c_7&=c_1(K) c_6(K)
     \cdot [ c_5(K,\Vert v_u\Vert, \Vert v_w\Vert)^2d(u,w)^2+\Vert v_u\Vert^2+\Vert v_w\Vert^2\Big].
\end{align*}
We use again Lemma \ref{lem:Kar2_a}, 
\begin{equation*}
        \Vert \Gamma_{u_+}^u \Exp_{u_+}^{-1}(w_+)) -\Exp_u^{-1}(w) -[v_u-\Gamma_w^u v_w] \Vert\leq (1+ c_3(K)R^2)\cdot c_7 d(u,w) .
\end{equation*}
Therefore we have linearized the iterate in $T_u\mM$. We should see how to transport it back to $T_{x}\mM$.
With Lemma \ref{lem:par_a} we have 
\begin{align*}
\Vert [\Gamma_u^{x}\Gamma_{u_+}^u -\Gamma_{u_+}^{x}] \Exp_{u_+}^{-1}(w_+)) \Vert =c_5(K)d(u,x)d(u_+,w_+)\Vert v_u\Vert.
\end{align*}
Note $v_u$ and $v_w$ are $-\eta\grad f(u)$ and $-\eta\grad f(w)$, we define $\nb v(x)$ the gradient of $v$, i.e., $-\eta H$. Using Hessian Lipschitz,
\begin{align*}
&\quad \Vert v_u-\Gamma_w^u v_w + \eta H(u) \Exp_u^{-1}(w)\Vert\\
 &=\Vert v_u-\Gamma_w^u v_w-\nabla v(u) \Exp_u^{-1}(w)\Vert \\
 &\leq \rho d(u,w)^2,
\end{align*}
and 
\begin{align*}
 &\quad \Vert \nabla v(u) \Exp_u^{-1}(w) -\Gamma_{x}^{u} \nabla v(x) \Gamma_u^{x}\Exp_u^{-1}(w)\Vert\leq \rho d(u,w)d(u,x).
\end{align*}
So we have 
\begin{align}
        &\quad \Vert \Gamma_{u_+}^x \Exp_{u_+}^{-1}(w_+) - (I+\nabla v(x)) \Gamma_u^{x}\Exp_u^{-1}(w) \Vert\notag \\
        &\leq c_7d(u,w)+ \rho d(u,w)(d(u,w)+d(u,x)) +c_5(K)d(u,x)d(u_+,w_+)\Vert v_u\Vert
        := D_1 \label{eq:u_w_compare}
\end{align}
\end{proof}

\subsection{Evolution of $\Expn_x(w) - \Expn_x(u)$}

We consider now the evolution of $\Expn_x(w) - \Expn_x(u)$ in the fixed tangent space $\mT_x\mM$. We show in the following lemma that it also approximately follows a linear iteration. 
\setcounter{lemma}{1}
\begin{lemma}\label{lemma_exp_inv_bound}
    Define $\gm = \sqrt{\hat\rho\ep}$, $\kappa = \frac{\bt}{\gm}$, and $\mathS = \sqrt{\eta\bt} \frac{\gm}{\hat\rho} \log^{-1}(\frac{d\ka}{\delta})$. 
Let us consider $x$ be a $(\ep,-\sqrt{\hat\rho\ep})$ saddle point, and  define $u^+ = \Exp_u(-\eta \grad f(u))$ and $w^+ = \Exp_w(-\eta \grad f(w))$.  Under Assumptions~\ref{assump:Lip_gra_a},~\ref{assump:Lip_Hes_a},~\ref{assump:sec_curv_a}, 
if all pairwise distances between $u,w,u^+,w^+,x$ are less than $12\mathS$, then for some explicit constant  $C(K,\rho,\bt)$ depending only on $K,\rho,\bt$, there is
\begin{align}
 &\quad \|\Expn_x(w^+) - \Expn_x(u^+) -(I-\eta H(x))(\Expn_x(w) - \Expn_x(u))\| \label{eq:exp_inv_bound}\\
 &\le C(K,\rho,\bt) d(u,w)\lc d(u,w)+d(u,x) + d(w,x)\rc. \notag
\end{align}
 \end{lemma}
This lemma controls the error of the linear approximation of the iterates hen mapped in $\mT_x\mM$ and largely follows from Lemma~\ref{lemma_2_pts_a}.
\begin{proof}
We have that
\begin{align}
 w &= \Exp_x(\Expn_x(w))\label{eq:w1}\\
 &= \Exp_u(\Expn_u(w)).\label{eq:w2}
\end{align}
Use \eq{w2}, let $a = \Expn_x(u)$ and $v=\Gm_u^x\Expn_u(w)$, Lemma \ref{lem:Kar1_a} suggests that
\begin{align*}
&\quad d(\Exp_{u}(\Expn_u(w)), \Exp_x(\Expn_x(u) + \Gm_u^x\Expn_u(w))) \\
&\le c_1(K)\|\Expn_u(w)\|(\|\Expn_u(w)\| + \|\Expn_x(u)\|)^2.
\end{align*}
Compare with \eq{w1}, we have 
\begin{align}
        &\quad d(\Exp_x(\Expn_x(w)), \Exp_x(\Expn_x(u) + \Gm_u^x\Expn_u(w)))\nonumber \\
&\le c_1(K)\|\Expn_u(w)\|(\|\Expn_u(w)\| + \|\Expn_x(u)\|)^2\nonumber \\
&:= D. \label{eq:def_D}
\end{align}
Denote the quantity above by $D$. Now use Lemma \ref{lem:Kar2_a}
\begin{align*}
   \|\Expn_x(w)- (\Expn_x(u) + \Gm_u^x\Expn_u(w)) \|\le (1+c_3(K)R^2)D.
\end{align*}
Analogously 
\begin{align*}
  \|\Expn_x(w_+)- (\Expn_x(u_+) + \Gm_{u_+}^x\Expn_{u_+}(w_+)) \|\le (1+c_3(K)R^2)D_+
\end{align*}
where 
\begin{equation}\label{eq:def_D+}
    \begin{split}
   D_+ =  c_1(K)\|\Expn_{u_+}(w_+)\|(\|\Expn_{u_+}(w_+)\|+ \|\Expn_x(u_+)\|)^2
\end{split}
\end{equation}
And we can compare $\Gm_{u}^x\Expn_{u}(w)$ and $\Gm_{u_+}^x\Expn_{u_+}(w_+)$ using \eq{u_w_compare}.
In the end we have 
\begin{equation*}
    \begin{split}
        &\quad\|\Expn_x(w^+) - \Expn_x(u^+)
        - (I-\eta H(x))(\Expn_x(w) - \Expn_x(u))\|\\
        &\le \|\Expn_x(w_+)- (\Expn_x(u_+) + \Gm_{u_+}^x\Expn_{u_+}(w_+)) \|\\
        &\quad + \|\Expn_x(w)- (\Expn_x(u) + \Gm_u^x\Expn_u(w)) \|\\
        &\quad + \|\Gm_{u_+}^x\Expn_{u_+}(w_+) - \Gm_u^x\Expn_u(w) - \nb v(x)\Gm_u^x\Expn_u(w) \|\\
        &\quad + \|\nb v(x)(\Gm_u^x\Expn_u(w) - (\Expn_x(w) - \Expn_x(u)))\|\\
        &\le (1+c_3(K)R^2)(D_++D)+ D_1 + \eta\|H(x)\|D.
    \end{split}
\end{equation*}
$D$, $D_+$ and $D_1$ are defined in \eq{def_D}, \eq{def_D+} and \eq{u_w_compare}, they are all order 
$
d(u,w)\big( d(u,w)+d(u,x) + d(w,x)\big) 
$
so we get the correct order in \eq{exp_inv_bound}.
  \end{proof}
\begin{figure}[t!]
      \centering
        \includegraphics[width = .5\textwidth]{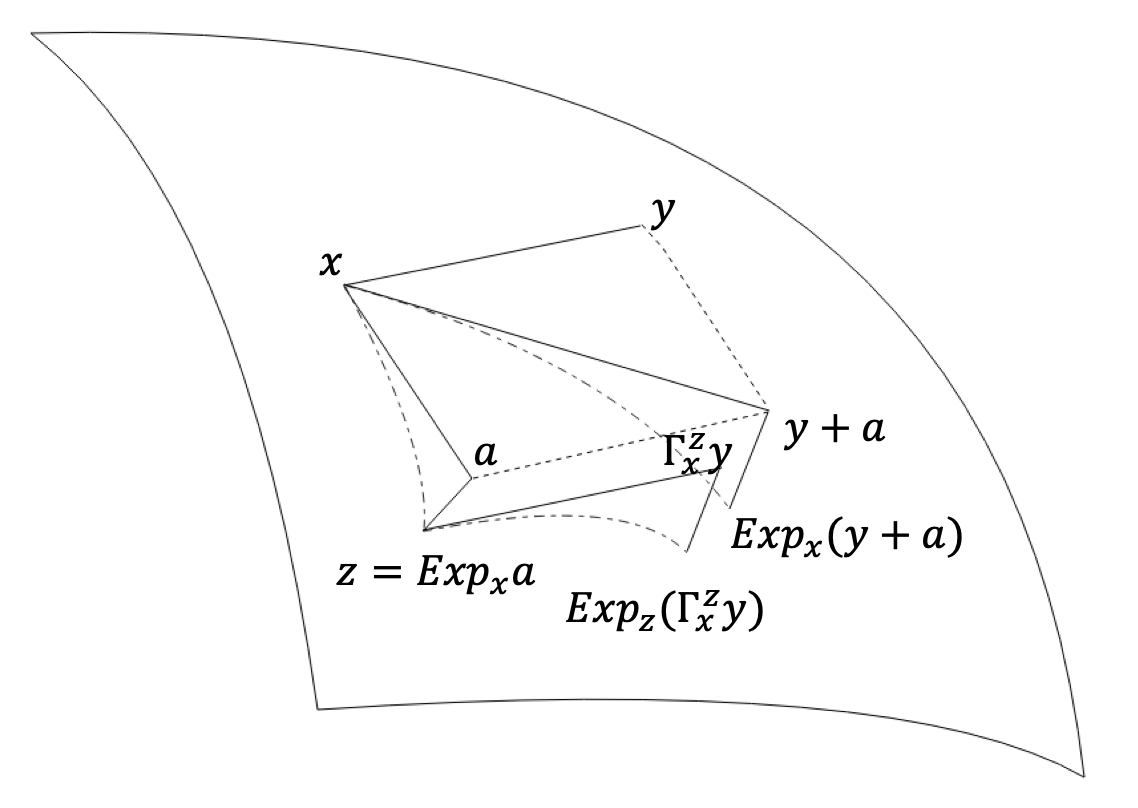}
      \caption{ Lemma \ref{lem:Kar1_a} bounds the difference of two steps starting from $x$: (1) take $y+a$ step in $\mT_x\mM$ and map it to manifold, and (2) take $a$ step in $\mT_x\mM$, map to manifold, call it $z$, and take $\Gm_x^z y$ step in $\mT_x\mM$, and map to manifold. $\Exp_z(\Gm_x^z y)$ is close to $\Exp_x(y+a)$.}\label{fig2}\vspace{-0.2cm}
    \end{figure}
\subsection{Control of two-steps iteration}
In the following lemma we control the distance between the point obtained after moving along the sum of two vectors in the tangent space, and the point obtained after moving a first time along the first vector and then a second time along the transport of the second vector. This is illustrated in Fig.~\ref{fig2}.
\begin{lemma}\label{lem:Kar1_a} 
 Let $x\in \mM$ and $y,a\in T_x\mM$. Let us denote by $z=\Exp_x(a)$ then under Assumption~\ref{assump:sec_curv_a}
\begin{equation}
\begin{split}
d(\Exp_x(y+a), \Exp_{z}(\Gamma_x^{z}y))\leq c_1(K) \min\{\Vert a \Vert, \Vert y \Vert\} (\Vert a \Vert+ \Vert y \Vert)^2. 
\end{split}
\end{equation}
\end{lemma}
This lemma which is crucial in the proofs of Lemma~\ref{lemma_exp_inv_bound} and Lemma~\ref{lemma_2_pts_a} tightens the result of \citet[][C2.3]{karcher}, which only shows an upper-bound 
$O(\Vert a \Vert(\Vert a \Vert+ \Vert y \Vert)^2)$.

\begin{proof}
We adapt the proof of \citet[][Eq. (C2.3) in App C2.2]{karcher}, the only difference being that we bound more carefully the initial normal component.  We restate here the whole proof for completeness. 

Let $x\in \mM$ and $y,a\in T_x\mM$. We denote by $\gamma(t)=\Exp_x(ta)$. We want to compare the point $\Exp_x(r (y+a))$ and $\Exp_\gamma(1)(\Gamma_x^{\gamma(1) y})$. These two points , for a fixed $r$ are joined by the curve 
\[
t\mapsto c(r,t) = \Exp_{\gamma(t)}(r\Gamma_x^{\gamma(t)}(y+(1-t)a)).
\]
We note that $\frac{d}{dt}c(r,t)$ is a Jacobi field along the geodesic $r\mapsto c(r,t)$, which we denote by $J_t(r)$. We importantly remark that the length of the geodesic $r\mapsto c(r,t)$ is bounded as $\|\frac{d}{dr}c(r,t)\|\leq \| y+(1-t)a\|$. We denote this quantity by $\rho_t= \| y+(1-t)a\|$. The initial condition of the Jacobi field $J_t$ are given by:
\begin{align*}
    &J_t(0)= \frac{d}{dt}\gamma(t)=\Gamma_x^{\gamma(t)}a\\
    &\frac{D}{dr}J_t(0)=\frac{D}{dr}\Gamma_x^{\gamma(t)}(y+(1-t)a)=-\Gamma_x^{\gamma(t)}a.
\end{align*}
These two vectors are linearly dependent and it is therefore possible to apply \citet[][Proposition A6]{karcher} to bound $J_t^{\rm norm}$. Moreover, following \citet[][App A0.3 ]{karcher}, the tangential component of the Jacobi field is known explicitly, independent of the metric, by
\[
J_t^{\rm tan}(r)=\left(J^{\rm tan}_t(0)+r\frac{D}{dr}J^{\rm tan}_t(0)\right) \frac{d}{dr}c(r,t)\]
where the initial conditions of the tangential component of the Jacobi fields are given by $J^{\rm tan}_t(0)=\langle J_t(0), \frac{\frac{d}{dr}c(r,t)}{\|\frac{d}{dr}c(r,t)\|}\rangle  $ and $\frac{D}{dr}J^{\rm tan}_t(0)=\langle \frac{D}{dr} J_t(0), \frac{\frac{d}{dr}c(r,t)}{\|\frac{d}{dr}c(r,t)\|}\rangle = -J^{\rm tan}_t(0) $. Therefore
\[
J_t^{\rm tan}(r)= (1-r) J^{\rm tan}_t(0) \frac{d}{dr}c(r,t),\]
and $J_t^{\rm tan}(1)=0$.

We estimate now the distance $d(\Exp_x(y+a), \Exp_{z}(\Gamma_x^{z}y))$  by the length of the curve $ t\mapsto c(r,t)$ as follows:
\begin{align*}
    d(\Exp_x(y+a), \Exp_{z}(\Gamma_x^{z}y)) \leq \int_0^1 \| \frac{d}{dt}c(1,t) \| dt =  \int_0^1 \| J_t^{\rm norm}(1) \| dt,
\end{align*}
where we use crucially that  $J_t^{tan}(1)=0$. 

We utilize \citep[][Proposition A.6]{karcher} to bound $\| J_t^{\rm norm}(1) \| $ as
\[
  \|J_t^{\rm norm}(1) \| \leq  \|J_t^{\rm norm}(0) \| (\cosh (\sqrt{K}\rho_t) -\frac{\sinh(\sqrt{K}\rho_t)}{\sqrt{K}\rho_t})
\]
using \citep[][Equation (A6.3)]{karcher} with $\kappa=0$, $f_\kappa(1)=0$ and recalling that the geodesics $r\mapsto c(r,t)$ have length $\rho_t$.
\begin{figure}[htbp!]
      \centering
        \includegraphics[width = .4\textwidth]{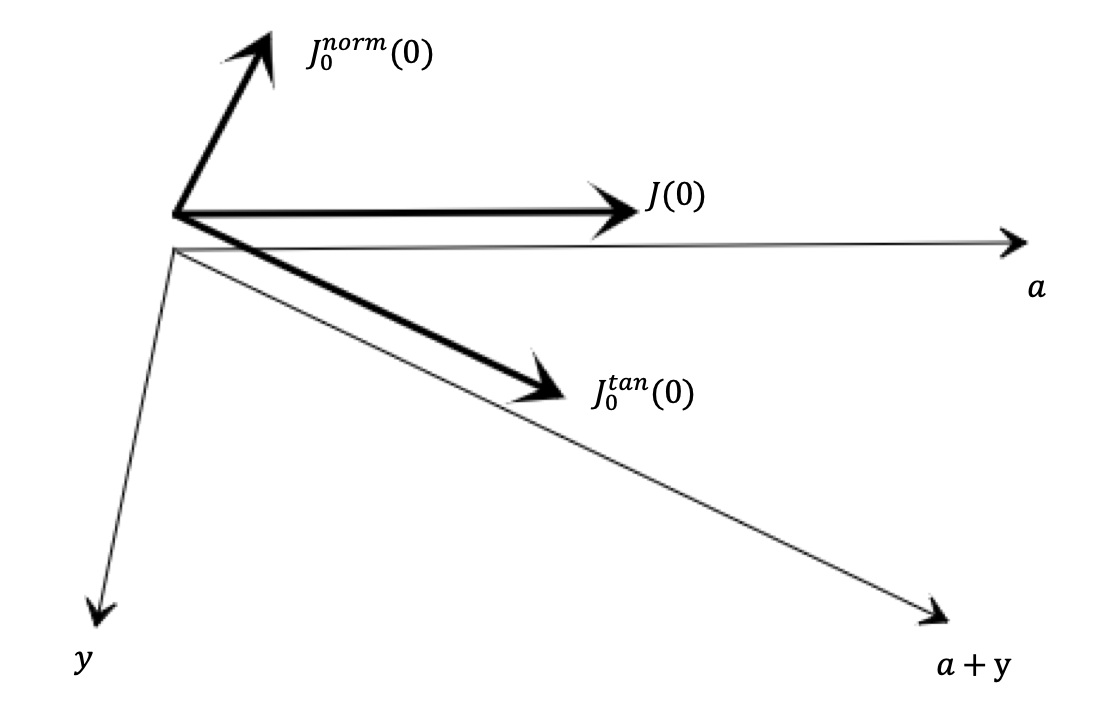}
      \caption{Figure for Lemma \ref{lem:Kar1_a}.}\label{fig3}
    \end{figure}

In particular for small value $\|a\|+\|y\|$ we have for some constant $c_1(K)$,
\[
  \|J_t^{\rm norm}(1) \| \leq  \|J_t^{\rm norm}(0) \| c_1(K) \rho_t^2.
\]

We  bound  $\|J_t^{\rm norm}(0) \|$ now. This is the main difference with the original proof of \citet{karcher} who directly bounded $\|J_t^{\rm norm}(0) \|\leq \|J_t(0) \|=\|a\|$ and $\rho_t\leq\|a\|+\|y\|$. Therefore his proof does not lead to the correct dependence in $\|y\|$. 

We have $J_t^0= \Gamma_x^{\gamma(t)}a$, and the tangential component (velocity of $r\ra c(r,t)$) is in the $\Gamma_x^{\gamma(t)}(y + (1-t)a)$ direction. Let $\tilde z = \Gamma_x^{\gamma(t)}(y + (1-t)a)$ and $\mP_{\tilde z^\perp}$ and $\mP_{a^\perp}$ denote the projection onto orthogonal complement of $\tilde z$ and $a$. 
\begin{align*}
    \|J_t^{\rm norm}(0)\|^2 &= \|\mP_{\tilde z^\perp}(a)\|^2\\
    &= \|a\|^2 - \frac{(a^T\tilde z)^2}{\|\tilde z\|^2}\\
    &= \frac{\|a\|^2}{\|\tilde z\|^2} \left(\|\tilde z\|^2 - \frac{(a^T\tilde z)^2}{\|\tilde z\|^2} \right)\\
    &\le \frac{\|a\|^2}{\|\tilde z\|^2} \|\mP_{a^\perp} (\Gamma_x^{\gamma(t)}(y+(1-t)a))\|^2\\
    &\le \frac{\|a\|^2}{\|\tilde z\|^2} \|\mP_{a^\perp} (\Gamma_x^{\gamma(t)}((1-t)a)) + \mP_{a^\perp} (\Gamma_x^{\gamma(t)} y)\|^2\\
    &= \frac{\|a\|^2}{\|\tilde z\|^2} \|\mP_{a^\perp} (\Gamma_x^{\gamma(t)} y)\|^2\\
    &\le \frac{\|a\|^2 \|y\|^2}{\|\tilde z\|^2}.
\end{align*}
So 
\begin{align*}
    \|J_t^{\rm norm}(1) \| &\le  \|J_t^{\rm norm}(0) \| c_1(K) \rho_t^2 \\
    &\le \frac{\|a\| \cdot \|y\|}{\|\tilde z\|} c_1(K)\|\tilde z\|^2 \\
    &\le c_1(K) \|a\| \cdot \|y\|(\|a\| + \|y\|),
\end{align*}
and 
\[
  d(\Exp_x(y+a), \Exp_{z}(\Gamma_x^{z}y)) \leq c_1(K) \|a\| \cdot \|y\|(\|a\| + \|y\|).
\]
\end{proof}

  \subsection{Auxilliary lemmas}

In the proofs of Lemma~\ref{lemma_2_pts_a} and Lemma~\ref{lemma_exp_inv_bound} we needed numerous auxiliary lemmas we are stating here.

We needed the following lemma which shows that both the exponential map  and its inverse are Lipschitz.
\begin{lemma}\label{lem:Kar2_a} Let $x,y,z\in M$, and the distance of each two points is no bigger than $R$. Then under Assumption~\ref{assump:sec_curv_a}
\begin{equation*}
    (1+c_2(K)R^2)^{-1}d(y,z)\leq \Vert \Exp_x^{-1}(y)-\Exp_x^{-1}(z)\Vert \leq(1+c_3(K)R^2)d(y,z).
\end{equation*}
\end{lemma}  
Intuitively this lemma relates the norm of the difference of two vectors of $\mT_x\mM$ to the distance between the corresponding points on the manifold $\mM$ and follows from bounds on the Hessian of the square-distance function~\citep[][Ex. 4 p. 154]{Sak96}.
\begin{proof}
The upper-bound is directly proven in \citet[][Proof of Cor. 1.6]{karcher}, and we prove the lower-bound via Lemma~\ref{lem:Kar1_a} in the supplement.
 Let $b = \Exp_y(\Gm_x^y(\Exp_x^{-1}(z)-\Exp_x^{-1}(y)))$. Using $d(y,b) = \|\Exp_y^{-1}(b)\|$ and Lemma \ref{lem:Kar1_a},
\begin{subequations}
\begin{align*}
    d(y,z) &\le d(y,b) + d(b, \Exp_x(\Expn_x(z)))\\
&\le \Vert \Exp_x^{-1}(y)-\Exp_x^{-1}(z)\Vert \\
&\quad + c_1(K) \Vert \Exp_x^{-1}(y)-\Exp_x^{-1}(z)\Vert (\Vert \Exp_x^{-1}(y)-\Exp_x^{-1}(z)\Vert+ \Vert \Exp_x^{-1}(y)\Vert)^2
\end{align*}
\end{subequations}
\end{proof}

The following contraction result is fairly classical and is proven using the Rauch comparison theorem from differential geometry~\citep{MR2394158}.
\begin{lemma}\label{lem:dis_a}
\citep[][Lemma 1]{rapidmix} Under Assumption~\ref{assump:sec_curv_a}, for $x,y\in \mM$ and $w\in T_x\mM$,
\[
d(\Exp_x(w),\Exp_y(\Gamma_x^y w))\leq c_4(K) d(x,y).
\]
\end{lemma}
Eventually we need the following corollary of the famous Ambrose-Singer holonomy theorem~\citep{ambrosesinger}.
\begin{lemma}\label{lem:par_a}
\citep[][Section 6]{karcher} Under Assumption~\ref{assump:sec_curv_a}, for $x,y,z\in \mM$ and $w\in T_x\mM$,
\[
\Vert \Gamma_y^{z}\Gamma_{x}^yw -\Gamma_{x}^{z}w \Vert \leq c_5(K) d(x,y)d(y,z) \Vert w\Vert.
\]
\end{lemma}

\section{Proof of Lemma \ref{lemma1_a} and \ref{lemma2_a}}\label{app:proof}
In this section we prove two important lemmas from which the proof of our main result mainly comes out. Then we show, in the last subsection, how to combine them to prove this main result.
\setcounter{lemma}{6}
\begin{lemma}\label{lemma1_a}
Assume Assumptions \ref{assump:Lip_gra_a}, \ref{assump:Lip_Hes_a}, \ref{assump:sec_curv_a}  hold, and 
\begin{equation}\label{eq:eps_bd2_a}
        \ep\le\min\left\{ \frac{\hat\rho}{56\max\{c_2(K),c_3(K)\}\eta\bt}\log\left(\frac{d\bt}{\sqrt{\hat\rho\ep}\delta}\right), \left(\frac{\mathI\hat\rho}{12\hat c\sqrt{\eta\bt}}\log\left(\frac{d\bt}{\sqrt{\hat\rho\ep}\delta}\right)\right)^2\right\}
\end{equation}
from the main theorem. There exists a constant $c_{\max}$, $\forall \hat c>3, \dt\in (0,\frac{d\ka}{e}]$, for any $u_{0}$ with $d(\tilde x, u_0)\le2\mathS/(\ka \log(\frac{d\ka}{\dt}))$, $\ka = \bt/\gm$.
\begin{equation*}
  T = \min\lbig \inf_t\lbig t|\tilde f_{u_0}(u_t) - f(u_0) \le -3\mF\rbig, \hat c \mathT  \rbig,
  \end{equation*}
  then $\forall \eta\le c_{\max}/\bt$, we have $\forall 0<t<T$, $d(u_0, u_t)\le 3(\hat c\mathS)$.
\end{lemma}
\begin{lemma}\label{lemma2_a}
Assume Assumptions \ref{assump:Lip_gra_a}, \ref{assump:Lip_Hes_a}, \ref{assump:sec_curv_a} and \eq{eps_bd2_a} hold. Take two points $u_0$ and $w_0$ which are perturbed from approximate saddle point, where $d(\tilde x, u_0)\le 2\mathS/(\ka \log(\frac{d\ka}{\dt}))$, 
$\Expn_{\tilde x}(w_0) - \Expn_{\tilde x}(u_0) = \mu r e_1$,
$e_1$ is the smallest eigenvector\footnote{``smallest eigenvector'' means the eigenvector corresponding to the smallest eigenvalue.} of $H(\tilde x)$, $\mu\in[\dt/(2\sqrt{d}),1]$, and the algorithm runs two sequences $\{u_t\}$ and $\{w_t\}$ starting from $u_0$ and $w_0$. Denote
\begin{equation*}
  T = \min\lbig \inf_t\lbig t|\tilde f_{w_0}(w_t) - f(w_0) \le -3\mF\rbig, \hat c \mathT  \rbig,
\end{equation*}
  then $\forall \eta\le c_{\max}/l$, if $\forall 0<t<T$, $d(\tilde x, u_t)\le 3(\hat c\mathS)$, we have $T<\hat c \mathT$.
\end{lemma}

\subsection{Proof of Lemma \ref{lemma1_a}}
Suppose $f(u_{t+1}) - f(u_t) \le -\frac{\eta}{2} \| \grad f(u_t)\|^2$.
\begin{subequations}
\begin{align*}
d (u_{\hat c \mathT}, u_0)^2 &\le (\sum\limits_{0}^{{\hat c \mathT}-1} d( u_{t+1} , u_t) )^2\\
&\le {\hat c \mathT} \sum\limits_{0}^{{\hat c \mathT}-1} d( u_{t+1}, u_t)^2\\
&\le \eta^2{\hat c \mathT}\sum\limits_{0}^{{\hat c \mathT}-1}\|\grad f(u_t)\|^2\\
&\le 2\eta{\hat c \mathT} \sum\limits_{0}^{{\hat c \mathT}-1} f(u_t) - f(u_{t+1})\\
&= 2\eta{\hat c \mathT} (f(u_0) - f(u_{\hat c \mathT}))\\
&\le 6\eta{\hat c \mathT}\mF = 6\hat c\mathS^2.
\end{align*}
\end{subequations} 
\subsection{Proof of Lemma \ref{lemma2_a}}
Note that, for any points inside a region with diameter $R$, under the assumption of Lemma \ref{lemma2_a}, we have $\max\{c_2(K),c_3(K)\}R^2\le 1/2$. 

Define 
$v_t = \Expn_{\tilde x}(w_t) - \Expn_{\tilde x}(u_t)$,
let $v_0 = e_1$ be the smallest eigenvector of $H(\tilde x)$, then let $\hat y_{2,t}$ be a unit vector, we have 
\begin{equation}\label{eq:v_iter}
\begin{split}
v_{t+1} &= (I - \eta H(\tilde x))v_t  + C(K,\rho,\bt) d(u_t,w_t)\\
&\quad\quad \cdot(d(u_t, \tilde x)+d(w_t, \tilde x) + d(\tilde x, u_0))\hat y_{2,t}.
\end{split}
\end{equation}
Let $C:=C(K,\rho,\bt)$. Suppose lemma \ref{lemma2_a} is false, then $0\le t\le T$, $d(u_t, \tilde x)\le 3\hat c\mathS$, $d(w_t, \tilde x)\le 3\hat c\mathS$, so $d(u_t, w_t)\le 6\hat c \mathS$, and the norm of the last term in \eq{v_iter} is smaller than $14\eta C\hat c\mathS \|v_t\|$. 

Lemma 4 in the main paper indicates that
\begin{equation}\label{eq:bound_on_v}
    \|v_t\|\in [1/2,2]\cdot d(u_t, w_t)= [3/2, 6]\cdot \hat c \mathS.
\end{equation}

Let $\psi_t$ be the norm of $v_t$ projected onto $e_1$, the smallest eigenvector of $H(0)$, and $\varphi_t$ be the norm of $v_t$ projected onto the remaining subspace. Then \eq{v_iter} is 
\begin{subequations}
\begin{align*}
\psi_{t+1}&\ge (1+\eta\gm)\psi_t - \mu \sqrt{\psi_t^2 + \phi_t^2},\\
\phi_{t+1}&\le (1+\eta\gm)\phi_t + \mu \sqrt{\psi_t^2 + \phi_t^2}.
\end{align*}
\end{subequations} 
Prove that for all $t\le T$, $\phi_t\le 4\mu t\psi_t$. Assume it is true for $t$, we have 
\begin{align*}
&4\mu(t+1) \psi_{t+1}
\ge 4\mu(t+1)\cdot \lc (1+\eta\gm)\psi_t - \mu \sqrt{\psi_t^2 + \phi_t^2}\rc,\\
&\phi_{t+1}\le 4\mu t(1+\eta\gm)\phi_t + \mu \sqrt{\psi_t^2 + \phi_t^2}.
\end{align*}
So we only need to show that 
\begin{equation*}
(1+4\mu(t+1))\sqrt{\psi_t^2 + \phi_t^2}\le (1+\eta\gm)\psi_t.
\end{equation*}
  
By choosing $\sqrt{c_{\max}}\le \frac{1}{56
\hat c^2}$ and $\eta \le c_{\max}/\bt$, we have 
\begin{equation*}
4\mu(t+1) \le 4\mu T \le 4\eta C \mathS \cdot 14\hat c^2 \mathT = 56\hat c^2\frac{C}{\hat\rho} \sqrt{\eta \bt} \le 1.
\end{equation*} 
This gives 
\begin{equation*}
4(1+\eta\gm)\psi_t \ge 2\sqrt{2\psi_t^2} \ge (1+4\mu(t+1))\sqrt{\psi_t^2 + \phi_t^2}.
\end{equation*} 
Now we know $\phi_t \le 4\mu t\psi_t\le \psi_t$, so
$
\psi_{t+1}\ge (1+\eta\gm)\psi_t -\sqrt{2}\mu\psi_t,
$
and 
\begin{equation*}
\mu = 14\hat c\eta C\mathS \le 14\hat c\sqrt{c_{\max}}\eta\gm C \log^{-1}(\frac{d\ka}{\dt})/\hat\rho\le \eta\gm/2,
\end{equation*}
so $\psi_{t+1}\ge(1+\eta\gm/2)\psi_t$.

We also know that $\|v_t\|\le 6\hat c\mathS$ for all $t\le T$ from \eq{bound_on_v}, so 
\begin{subequations}
\begin{align*}
6\hat c\mathS&\ge \|v_t\|\ge\psi_{t}\ge (1+\eta\gm/2)^t\psi_0\\
&=(1+\eta\gm/2)^t\frac{\mathS}{\ka}\log^{-1}(\frac{d\ka}{\dt})\\
&\ge (1+\eta\gm/2)^t\frac{\dt\mathS}{2\sqrt{d}\ka}\log^{-1}(\frac{d\ka}{\dt}).
\end{align*}
\end{subequations}
This implies 
\begin{subequations}
\begin{align*}
T&< \frac{\log(12\frac{\ka\sqrt{d}}{\dt} \hat c \log(\frac{d\ka}{\dt}))}{2\log (1+\eta\gm/2)}\\
&\le  \frac{\log(12\frac{\ka\sqrt{d}}{\dt} \hat c \log(\frac{d\ka}{\dt}))}{\eta\gm}\\
&\le (2+\log(12\hat c))\mathT.
\end{align*}
\end{subequations}
By choosing $\hat c$ such that $2+\log(12\hat c)<\hat c$, we have $T\le \hat c\mathT$, which finishes the proof. 
\subsection{Proof of function value decrease at an approximate saddle point}
With Lemma \ref{lemma1_a} and \ref{lemma2_a} proved, we can lower bound the function value in $O(\mathT)$ iterations decrease by $\Omega(\mF)$, thus match the convergence rate in the main theorem. Let $T^\prime := \inf_t\lbig t|\tilde f_{u_0}(u_t) - f(u_0) \le -3\mF\rbig$. Let $\widecheck{\ }$ denote the operator $\Exp_{u_0}^{-1}(\cdot)$. If $T^\prime\le T$,
\begin{subequations}
\begin{align*}
&\quad f(u_{T^\prime}) - f(u_0) \\
&\le \nb f(u_0)^T(u_{T^\prime} - u_0) + \half H(u_0)[\wc u_{T^\prime} - u_0, \wc u_{T^\prime} - u_0] \\
&\quad  +\frac{\rho}{6} \|\wc u_{T^\prime} - u_0\|^3\\
&\le \tilde f_{u_0}(u_t) - f(u_0) + \frac{\rho}{2}d(u_0, \tilde x) \|\wc u_{T^\prime} - u_0\|^2\\ 
&\le -3\mF + O(\rho\mathS^3) \le -2.5\mF.
\end{align*}
\end{subequations} 
If $T^\prime>T$, then $\inf_t\lbig t|\tilde f_{w_0}(w_t) - f(w_0) \le -3\mF\rbig\le T$, and we know $f(w_T)-f(w_0)\le -2.5\mF$.
\begin{remark}
What is left is bounding the volume of the stuck region, to get the probability of getting out of the stuck region by the perturbation. The procedure is the same as in \citet{GeSham}. We sample from a unit ball in $\mathcal{T}_x\mM$, where $x$ is the approximate saddle point. In Lemma \ref{lemma1_a} and \ref{lemma2_a}, we study the inverse exponential map at the approximate saddle point $x$, and the coupling difference between $\Expn_x(w)$ and $\Expn_x(u)$. The iterates we study and the noise are all in the tangent space $\mathcal{T}_x\mM$ which is a Euclidean space, so the probability bound is same as the one in \citet{GeSham}.
\end{remark}


\end{document}